\documentclass[a4paper,8pt]{article}

\usepackage{hyperref}

\usepackage{amsmath,amsfonts,amssymb,amsthm,graphicx,psfrag,subfigure,mathtools,dsfont,enumerate}
\usepackage[german,english]{babel}
\usepackage[utf8]{inputenc}	
\usepackage{xparse}
\usepackage[backend=bibtex,style=numeric,maxnames=5]{biblatex}
\newtheorem{satz}{Satz}[section]
\newtheorem{theorem}[satz]{Theorem}
\newtheorem{lemma}[satz]{Lemma}

\theoremstyle{definition}
\newtheorem{definition}[satz]{Definition}

\newtheorem{bemerkung}[satz]{Remark}
\newtheorem{annahme}[satz]{Assumptions}
{\begin{proof}[Proof]}
	{\end{proof}}

\makeatletter                    
\def\blfootnote{\xdef\@thefnmark{}\@footnotetext}
\DeclareMathAlphabet\mathbfcal{OMS}{cmsy}{b}{n}

\usepackage{xparse}

\makeatletter

\DeclarePairedDelimiter{\norm}{\lVert}{\rVert}
\NewDocumentCommand{\norml}{ s O{} m m }{%
	\IfBooleanTF{#1}{\norm*{#4}}{\norm[#2]{#4}}_{L^{#3}}%
}
\NewDocumentCommand{\normL}{ s O{} m m }{%
	\IfBooleanTF{#1}{\norm*{#4}}{\norm[#2]{#4}}_{\L^{#3}}%
}
\NewDocumentCommand{\normw}{ s O{} m m m }{%
	\IfBooleanTF{#1}{\norm*{#4}}{\norm[#2]{#5}}_{W^{#3,#4}}%
}
\NewDocumentCommand{\normW}{ s O{} m m m }{%
	\IfBooleanTF{#1}{\norm*{#4}}{\norm[#2]{#5}}_{\W^{#3,#4}}%
}
\NewDocumentCommand{\normh}{ s O{} m m }{
	\IfBooleanTF{#1}{\norm*{#4}}{\norm[#2]{#4}}_{H^{#3}}
}
\NewDocumentCommand{\normH}{ s O{} m m }{%
	\IfBooleanTF{#1}{\norm*{#4}}{\norm[#2]{#4}}_{\H^{#3}}%
}
\makeatother

\def\R{\mathbb R}

\def\N{\mathbb N}
\def\n{\mathbf n}
\def\a{\mathbf a}

\def\f{\mathbf{f}}

\def\L{\mathbf L}
\def\H{\mathbf H}

\def\u{\mathbf{u}}
\def\v{\mathbf{v}}

\def\T{\mathbf{T}}
\def\I{\mathbf{I}}
\def\D{{\mathbf{D}}}

\def\W{{\mathbf W}}

\def\q{{\mathbf{q}}}

\def\intO{\int_\Omega}
\def\intT{\int_{0}^{T}}

\def\dx{\;\mathrm dx}

\def\dt{\;\mathrm dt}

\def\del{\partial}

\def\delt{\partial_{t}}

\def\deln{\partial_\n}

\def\d{{\mathrm{d}}}
\DeclareMathOperator*{\esssup}{ess\,sup}

\def\grad{\nabla}
\def\laplace{\Delta}
\def\divergence{\textnormal{div}}

\bibliography{On_a_Cahn_Hilliard_Brinkman_model_for_tumour_growth_and_its_singular_limits}

\title{On a Cahn-Hilliard-Brinkman model for tumour growth and its singular limits}
\author{Matthias Ebenbeck\quad Harald Garcke}
\author{  
		Matthias Ebenbeck\footnote{Fakult\"at f\"ur Mathematik,  
		Universit\"at Regensburg,
		93040 Regensburg,
		Germany, e-mail: {\sf matthias.ebenbeck@ur.de}}\qquad
		Harald Garcke\footnote{Fakult\"at f\"ur Mathematik,  
		Universit\"at Regensburg,
		93040 Regensburg,
		Germany, e-mail: {\sf harald.garcke@mathematik.uni-regensburg.de}}
		}
\date{}

\begin{document}
\maketitle
\begin{abstract}
	In this work, we study a model consisting of a Cahn-Hilliard-type equation for the concentration of tumour cells coupled to a reaction-diffusion type equation for the nutrient density and a Brinkman-type equation for the velocity. We equip the system with Neumann boundary for the tumour cell variable and the chemical potential, Robin-type boundary conditions for the nutrient and a \grqq no-friction" boundary condition for the velocity, which allows us to consider solution dependent source terms. Well-posedness of the model as well as existence of strong solutions will be established for a broad class of potentials. We will show that in the singular limit of vanishing viscosities we recover a Darcy-type system related to Cahn-Hilliard-Darcy type models for tumour growth which have been studied earlier. An asymptotic limit will show that the results are also valid in the case of Dirichlet boundary conditions for the nutrient.
\end{abstract}

\noindent{\bf AMS-Classification:} 
35K35, 
35Q92, 
92C50, 
35D30, 
76D07 

\section{Introduction}
\numberwithin{equation}{section}
The growth of living cancer cells is affected by many biological and chemical mechanisms. Although the amount of experimental data coming from clinical experiments is quite extensive, the understanding of involved mechanisms and biological effects is still at an unsatisfying level. In the recent past, several mathematical models for tumour growth have been developed and simulated and some of them seem to compare well with tumour data coming from clinical experiments, see \cite{AgostiEtAl,BearerEtAl,FrieboesEtAl}. These models may provide further insights into tumour growth dynamics to understand key mechanisms and to design new treatment strategies. \newline
Many models are based on the hypothesis that different tissue components of the tumour (viable, quiescent, necrotic) are separated by a sharp interfacial layer and therefore can be described by free boundary problems, see \cite{ByrneChaplain2,FranksKing,Friedman,Friedman3}. As a young tumour does not have its own vascular system and must therefore consume growth factors like nutrients or oxygen from the surrounding host tissue, in the early stage of growth the tumour may undergo morphological instabilities like fingering or folding (see \cite{CristiniLowengrubNie,CristiniEtAl}) to grow without angiogenesis and to overcome diffusional limitations. This leads to highly challenging mathematical problems when modelling the tumour in the context of free boundary problems since changes in topology have to be tracked.\newline 
To overcome this difficulties, it has turned out that diffuse interface models, treating the tumour as a collection of cells, are a good strategy to describe the evolution and interactions of different species. These models are typically based on a multiphase approach, on balance laws for the single constituents, like mass and momentum balance, on constitutive laws and on thermodynamic principles. Several additional variables describing the extracellular matrix (ECM), growth factors or inhibitors, can be incorporated into these models, biological mechanisms like chemotaxis, apoptosis or necrosis and effects of stress, plasticity or viscoelasticity can be included, see \cite{AstaninPreziosi,GarckeLamNuernbergSitka,OdenTinsleyHawkins,PreziosiTosin}.  Many contributions in the literature consider a mixture of two components (healthy and surrounding tissue), modelled as a two-phase flow and coupled to a reaction-diffusion type equation for an unknown species acting as a nutrient for the tumour, like oxygen or glucose, leading to coupled systems of partial differential equations of Cahn-Hilliard-type, see \cite{ColliGilardiHilhorst,FrigeriGrasselliRocca,GarckeLam3,GarckeLam2,HawkinsZeeKristofferOdenTinsley,HilhorstKampmannNguyenZee}. Some of them incorporate effects based on the momentum balance equation which leads to models involving e.g. Darcy's law for the velocity or a Stokes like equation, see \cite{EbenbeckGarcke,GarckeLam1,GarckeLam4,JiangWuZheng}.
\newline\newline
In the following, we will consider a Cahn-Hilliard-Brinkman system for tumour growth. For a bounded domain $\Omega\subset \R^d,\,d=2,3,$ and a fixed time $T>0$, we consider for $Q\coloneqq \Omega\times (0,T)$ the following system of equations
\begin{subequations}
	\label{MEQ}
	\begin{alignat}{3}
	\label{MEQ_1}\divergence(\v)&=\Gamma_{\v}(\varphi,\sigma)&& \quad \text{in }   Q,\\
	\label{MEQ_2}-\divergence(\T(\v,p)) +\nu\v&= (\mu + \chi\sigma)\grad\varphi  && \quad \text{in } Q,\\
	\label{MEQ_3}\delt\varphi + \divergence(\varphi\v) &= \divergence(m(\varphi)\grad\mu)+\Gamma_{\varphi}(\varphi,\sigma)&& \quad \text{in } Q,\\
	\label{MEQ_4}\mu&= \epsilon^{-1}\psi'(\varphi)-\epsilon\laplace\varphi -\chi\sigma && \quad \text{in } Q,\\
	\label{MEQ_5}0 &= \laplace\sigma - h(\varphi)\sigma&& \quad \text{in } Q,
	\end{alignat}
\end{subequations}
where the viscous stress tensor is defined by
\begin{equation}
\label{DEF_STRTENS}\T(\v,p) \coloneqq 2\eta(\varphi) \D\v+\lambda(\varphi)\divergence(\v)\I  - p\I ,
\end{equation}
and the symmetric velocity gradient is given by
\begin{equation*}
\D\v\coloneqq \frac{1}{2}(\grad\v+\grad\v^T).
\end{equation*}
In \eqref{MEQ}-\eqref{DEF_STRTENS}, $\v$ denotes the volume-averaged velocity of the mixture, $p$ denotes the pressure, $\sigma$ denotes the concentration of an unknown species acting as a nutrient, $\varphi\in[-1,1]$ denotes the difference in volume fractions,
with $\{\varphi = 1\}$ representing the unmixed tumour tissue and  $\{\varphi = -1\}$ representing the surrounding healthy tissue, and $\mu$ denotes the chemical potential for $\varphi$. The function $m(\cdot)$ represents the mobility for the phase variable $\varphi$. The constant $\epsilon >0$ is related to the thickness of
the diffuse interface, whereas $\nu$ is a positive constant representing the fluid permeability. Moreover, the functions $\eta(\cdot)$ and $\lambda(\cdot)$ are non-negative and represent the shear and the bulk viscosity, respectively. The chemotaxis parameter $\chi$ is a non-negative constant.\newline
By $\n $ we will denote the outer unit normal on $\del\Omega$, and $\deln g \coloneqq \grad g\cdot \n $ is the directional derivative .
We equip the system with the following initial and boundary conditions
\begin{subequations}
	\label{BIC}
	\begin{alignat}{3}
	\label{BC_1}\deln\mu=\deln\varphi &= 0 &&\quad \text{on }\del\Omega\times (0,T)\eqqcolon\Sigma,\\
	\label{BC_2}\deln\sigma &= K(\sigma_{\infty}-\sigma)&&\quad\text{on }\Sigma,\\
	\label{BC_3}\T(\v,p)\n  &= \mathbf{0}&&\quad\text{on }\Sigma,\\
	\label{IC}\varphi(0) &= \varphi_0 &&\quad\text{in }\Omega,
	\end{alignat} 
\end{subequations}
where $\varphi_0,\,\sigma_{\infty}$ are given functions and $K$ is a positive permeability constant.
We now motivate and compare our model with other models in the literature.
\begin{enumerate}[$\bullet$]
	\item In \cite{EbenbeckGarcke}, the authors of this paper considered a similar model with \eqref{MEQ_2}, \eqref{MEQ_5} and \eqref{BC_2} replaced by
	\begin{subequations}
		\begin{alignat}{3}
		\label{MEQ_NONSTAT_1}-\divergence(\T(\v,p))+\nu\v &= \mu\grad\varphi + (\chi_{\sigma} + \chi(1-\varphi))\grad\sigma&&\quad\text{in }Q,\\
		\label{MEQ_NONSTAT_2}\delt\sigma + \divergence(\sigma\v) &= \divergence(n(\varphi)(\chi_{\sigma}\grad\sigma-\chi\grad\varphi))-\Gamma_{\sigma}(\varphi,\sigma,\mu)&&\quad \text{in }Q,\\
		\label{MEQ_NONSTAT_3}n(\varphi)\chi_{\sigma}\deln\sigma &= K(\sigma_{\infty}-\sigma)&&\quad\text{on }\Sigma,
		\end{alignat}
	\end{subequations}
	where $\chi_{\sigma}$ is a given positive constant denoting the nutrient diffusion coefficient, $n(\varphi)$ is a mobility function and $\Gamma_{\sigma}$ is a source or sink term similar to ours.  Reformulating the pressure as
	\begin{equation*}
	\tilde{p}\coloneqq p - \frac{\chi_{\sigma}}{2}|\sigma|^2 + \chi\sigma(1-\varphi),
	\end{equation*}
	we can indeed recover \eqref{MEQ_2} from \eqref{MEQ_NONSTAT_1}. For a discussion regarding pressure reformulations, we refer to \cite{GarckeLam1}.
	Following the arguments in \cite[Sec. 3.3.3]{GarckeLamSitkaStyles}, we introduce an active transport parameter $\kappa$, a new mobility function $\mathcal{D}(\varphi)$ and set 
	\begin{equation*}
	n(\varphi) = \kappa D(\varphi)\chi^{-1},\quad \chi_{\sigma}= \kappa^{-1}\chi,
	\end{equation*}
	which allows us to \grqq decouple" active transport and chemotaxis mechanisms. Neglecting mechanism due to active transport, we can rewrite \eqref{MEQ_NONSTAT_2}-\eqref{MEQ_NONSTAT_3} as
	\begin{subequations}
		\begin{alignat}{3}
		\label{MEQ_NONSTAT_2a}\delt\sigma + \divergence(\sigma\v) &= \divergence(\mathcal{D}(\varphi)\grad\sigma)-\Gamma_{\sigma}(\varphi,\sigma,\mu)&&\quad\text{in }Q,\\
		\label{MEQ_NONSTAT_3a}\mathcal{D}(\varphi)\deln\sigma &= K(\sigma_{\infty}-\sigma)&&\quad\text{on }\Sigma.
		\end{alignat}
	\end{subequations}
	Performing a non-dimensionalisation and using that the ratio between nutrient diffusion time-scale (minutes or hours) and cell-proliferation time-scale (days or weeks) is quite small, we can drop the time derivative and the convection term in \eqref{MEQ_NONSTAT_2a} to obtain
	\begin{subequations}
		\begin{alignat}{3}
		\label{MEQ_NONSTAT_2b}0 &= \laplace\sigma-\Gamma_{\sigma}(\varphi,\sigma,\mu)&&\quad\text{in }Q,\\
		\label{MEQ_NONSTAT_3b}\deln\sigma &= K(\sigma_{\infty}-\sigma)&&\quad\text{on }\Sigma.
		\end{alignat}
	\end{subequations}
	We omit the details and refer to \cite{GarckeLam3,Preziosi,WiseLowengrubFrieboesChristini}.
	\item Brinkman's law was first proposed in \cite{Brinkman} and has been derived rigorously by several authors using a homogenization argument for the Stokes equation, see \cite{Allaire1,SanchezPalencia}. It can be interpreted as an interpolation between Stokes flow and Darcy's law since the former one is approximated on small length scales whereas the latter one on large length scales, see \cite{DurlofskyBrady}. Setting $\sigma = 0$, neglecting source terms, i.\,e. $\Gamma_{\varphi} = \Gamma_{\v}=0$, and imposing Dirichlet boundary conditions for the velocity, equations \eqref{MEQ}, \eqref{BIC} have been studied in \cite{BosciaContiGrasselli} for smooth potentials and in \cite{ContiGiorgini} for logarithmic potentials. However, if source terms are present, the Dirichlet boundary condition leads to a compatibility condition since
	\begin{equation*}
	\intO \Gamma_{\v}\dx  = \intO \divergence(\v)\dx  = \int_{\del\Omega}\v\cdot\n \d\mathcal{H}^{d-1} = 0,
	\end{equation*}
	which can in general not be fulfilled if $\Gamma_{\v}$ depends on additional variables like $\varphi$ and $\sigma$. The \grqq no-friction" boundary condition does not lead to such an condition and therefore enables us to consider solution dependent source terms. 
	\item In the case $\v = 0$ and with \eqref{MEQ_5} replaced by \eqref{MEQ_NONSTAT_2}, the model \eqref{MEQ}, \eqref{BIC} was analysed in \cite{GarckeLam3} and well-posedness has been shown in the case of Dirichlet boundary conditions for $\varphi,\mu,\sigma$ in \cite{GarckeLam2}, where they could rigorously prove the quasi-static limit. In \cite{FrigeriGrasselliRocca}, well-posedness has been established for $\deln\varphi = \deln\mu = \deln\sigma = 0$ on $\Sigma$, for a large class of potentials and with source terms of the form
	\begin{equation*}
	\Gamma_{\varphi}(\varphi,\mu,\sigma) = P(\varphi)(\sigma- \mu),\quad \Gamma_{\sigma}(\varphi,\mu,\sigma) = P(\varphi)(\sigma- \mu),
	\end{equation*}
	where $P(\cdot)$ is a nonlinear proliferation function. 
	\item When setting $\eta(\cdot)=\lambda(\cdot)\equiv 0$, the model reduces to the so called Cahn-Hilliard-Darcy model. In the absence of nutrients, this model has been studied in \cite{JiangWuZheng} for $\deln\varphi = \deln\mu = \v\cdot\n = 0$ on $\Sigma$ and with prescribed source terms $\Gamma_{\varphi} = \Gamma_{\v} = \mathcal{S}$. Setting $\sigma = 0$, prescribing $\deln\varphi = \deln\mu = \v = 0$ on $\Sigma$ and setting $\Gamma_{\v}=\Gamma_{\varphi}=0$, solutions of the Cahn-Hilliard-Darcy system have been established upon considering the zero viscosity limit in the Brinkman equation, see \cite{BosciaContiGrasselli}. Furthermore, we want to refer to \cite{DaiEtAl}, where they considered a multispecies Cahn-Hilliard-Darcy model for tumour growth with quasi-static nutrient equation.
\end{enumerate}
In the following, we want to outline some challenges arising in the analysis. When testing the Brinkman equation with $\v$, we have to estimate the term $\intO p\divergence(\v)\dx$. Hence, we need to get an estimate on $\norml{2}{p}$ in the absence of any a-priori-estimates. To overcome this difficulty, we will use the so-called method of \grqq Subtracting the divergence". More precisely, we choose $\v-\u $ as a test-function in the momentum equation, where $\u $ satisfies
\begin{equation*}
\divergence(\u) = \Gamma_{\v}\quad\text{in }\Omega,\quad \u = \frac{1}{|\del\Omega|}\left(\intO \Gamma_{\v}\right)\n \quad\text{on }\del\Omega.
\end{equation*}
Although this prevents us from estimating the pressure, we now have to handle the term
\begin{align*}
\intO\mu\grad\varphi\cdot\u \dx &= \intO (\mu-\mu_{\Omega})\grad\varphi\cdot\u \dx + \mu_{\Omega}\intO \grad\varphi\cdot\u\dx \\
&= \intO (\mu-\mu_{\Omega})\grad\varphi\cdot\u \dx + \mu_{\Omega}\left(\frac{1}{|\del\Omega|}\intO \Gamma_{\v}\dx\int_{\del\Omega}\varphi\d \mathcal{H}^{d-1} - \intO \varphi\Gamma_{\v}\dx\right),
\end{align*}
where $\mu_{\Omega} = \frac{1}{|\Omega|}\intO \mu\dx$. To handle the boundary integral, we will derive an estimate for the $L^{\rho}(\del\Omega)-$norm for $\varphi$, where $\rho\in [2,6]$ is an exponent connected to the growth rate of the potential $\psi(\cdot)$.\newline
Furthermore, we comment on the assumption $\sigma_{\infty}\in L^4(L^2(\del\Omega))$, which is not needed to prove existence of weak solutions, but crucial to establish well-posedness of the system. Indeed, this enables us to estimate the velocity in $L^{\frac{8}{3}}(0,T;\H^1)$ (see Section 3.3). As a consequence, we can handle the term
\begin{equation*}
\intO 2(\eta(\varphi_1)-\eta(\varphi_2))\D\v_2\colon\grad\v\dx 
\end{equation*}
in the proof for continuous dependence (see Section 5). We remark that this term does not arise in the case of constant viscosity.\newline
Finally, in the proof for existence of strong solutions we will derive an estimate for the time derivative of the nutrient concentration by using a difference quotient method. This is needed due to the fact that the $L^2$-orthogonal projection $\mathbb{P}_n$ onto the $n$-dimensional Galerkin solution spaces is not continuous on the whole space $H^2$. Indeed, when testing \eqref{MEQ_4} with $\laplace\delt\varphi$ in the Galerkin scheme and integrating by parts twice, we encounter the term
\begin{equation*}
\intT\intO \chi\laplace(\mathbb{P}_n\sigma)\delt\varphi\d x\d t.
\end{equation*} 
Although we can control $\sigma\in L^2(H^2)$, an estimate of $\laplace\mathbb{P}_n\sigma\in L^2(L^2)$ can not be deduced due to \eqref{BC_2}. If the time derivative of $\sigma$ fulfils $\delt\sigma\in L^2(H^1)$, a control of $\laplace\mathbb{P}_n\sigma\in L^2(L^2)$ is not needed, see Section 6.
\subsection{Notation and preliminaries}
We first want to fix some notation: For a (real) Banach space $X$ we denote by $\norm{\cdot}_X$ its norm, by $X^*$ the dual space and by $\langle \cdot{,}\cdot \rangle_X$ the duality pairing between $X^*$ and $X$. For an inner product space $X$, the inner product is denoted by $(\cdot{,}\cdot)_X$. We define the scalar product of two matrices by
\begin{equation*}
\mathbf{A}\colon \mathbf{B}\coloneqq \sum_{i=1}^{d}a_{jk}b_{jk}\quad\text{for } \mathbf{A},\mathbf{B}\in\R^{d\times d}.
\end{equation*}
For the standard Lebesgue and Sobolev spaces with $1\leq p\leq \infty$, $k>0$,  we use the notation $L^p\coloneqq L^p(\Omega)$ and $W^{k,p}\coloneqq W^{k,p}(\Omega)$ with norms $\norml{p}{\cdot}$ and $\normw{k}{p}{\cdot}$ respectively.  In the case $p=2$ we use $H^k\coloneqq W^{k,2}$ and the norm $\normh{k}{\cdot}$. For $\beta\in (0,1)$ and $r\in (1,\infty)$, we will denote the Lebesgue and Sobolev spaces on the boundary by $L^p(\del\Omega)$ and $W^{\beta,r}(\del\Omega)$ with corresponding norms $\norm{\cdot}_{L^p(\del\Omega)}$ and $\norm{\cdot}_{W^{\beta,r}(\del\Omega)}$ (see \cite[Chap. I.3]{Sohr} for more details). By $\L^p$, $\W^{k,p}$, $\H^k$, $\L^p(\del\Omega)$ and $\W^{\beta,r}(\del\Omega)$, we will denote the corresponding spaces of vector valued and matrix valued functions. We denote the space $H_0^1$ as the completion of $C_0^{\infty}(\Omega)$ with respect to the $H^1$-norm. For the Bochner spaces, we use the notation $L^p(X)\coloneqq L^p(0,T;X)$ for a Banach space $X$ with $p\in [1,\infty]$.  For the dual space $X^*$ of a Banach space $X$, we introduce the (generalised) mean value by 
\begin{equation*}
v_{\Omega}\coloneqq \frac{1}{|\Omega|}\intO v\dx \quad\text{for } v\in L^1, \quad v_{\Omega}^*\coloneqq \frac{1}{|\Omega|}\langle v{,}1\rangle_X\quad\text{for } v\in X^*.
\end{equation*}
Moreover, we introduce the function spaces
\begin{align*}
&L_0^2\coloneqq \{w\in L^2\colon w_{\Omega}=0\},\quad H_N^2\coloneqq \{w\in H^2\colon \deln w = 0 \text{ on } \del\Omega\},\\
&(H^1)_0^*\coloneqq \{f\in (H^1)^*\colon f_{\Omega}^* =0\}.
\end{align*}
Then, the Neumann-Laplace operator $-\laplace_N\colon H^1\cap L_0^2\to (H^1)_0^*$ is positive definite and self-adjoint. In particular, by the Lax-Milgram theorem and Poincaré's inequality, the inverse operator $(-\laplace_N)^{-1}\colon (H^1)_0^*\to H^1\cap L_0^2$ is well-defined, and we set $u\coloneqq (-\laplace_N)^{-1}f$ for $f\in (H^1)_0^*$ if $u_{\Omega}=0$ and 
\begin{equation*}
-\laplace u = f\text{ in }\Omega,\quad\deln u =0\text{ on }\del\Omega.
\end{equation*}
We have dense and continuous embeddings $H_N^2\subset H^1\subset L^2\simeq (L^2)^*\subset (H^1)^*\subset (H_N^2)^*$ and the identifications $\langle u{,}v\rangle_{H^1}=(u{,}v)_{L^2}$, $\langle u{,}w\rangle _{H^2} = (u{,}w)_{L^2}$ for all $u\in L^2,\, v\in H^1$ and $w\in H_N^2$. We also want to state the following inequalities resulting from elliptic regularity theory and integration by parts:
\begin{subequations}\label{ELLIPTIC_EST}
	\begin{alignat}{3}
	\label{ELLIPTIC_EST_1}\norml{2}{\laplace f}&\leq \normL{2}{\grad f}^{\frac{1}{2}}\normL{2}{\grad\laplace f}^{\frac{1}{2}}\leq C\normh{1}{\varphi}^{\frac{1}{2}}\normh{3}{\varphi}^{\frac{1}{2}}&&\qquad\forall f\in H_N^2\cap H^3,\\
	\label{ELLIPTIC_EST_2}\normh{2}{f}&\leq C\left(\norml{2}{f}+\norml{2}{\laplace f}\right)&&\qquad\forall f\in H_N^2,
	\end{alignat}
\end{subequations}
with a constant $C$ depending only on $\Omega$.
Furthermore, we define
\begin{equation*}
L_{\divergence}^2(\Omega)\coloneqq \left\{ \f\in \L^2\colon \divergence(\f)\in L^2 \right\}
\end{equation*}
equipped with the norm
\begin{equation*}
\norm{\f}_{L_{\divergence}^2(\Omega)}\coloneqq \left(\normL{2}{\f}^2 + \norml{2}{\divergence(\f)}^2\right)^{\frac{1}{2}},
\end{equation*}
which is a reflexive Banach space (see \cite{SimaderSohr}). In particular, for $\u\in L_{\divergence}^2(\Omega)$ we have
\begin{equation}
\label{TRACE_GENERAL}\langle \u\cdot\n{,}\Phi\rangle_{H^{\frac{1}{2}}(\del\Omega)} = \intO\u\cdot\nabla\Phi\d x + \intO \Phi\divergence(\u)\d x\quad \forall \Phi\in H^1
\end{equation}
and 
\begin{equation}
\label{TRACE_ESTIM}\norm{\u\cdot\n}_{\left(H^{\frac{1}{2}}(\del\Omega)\right)^*}\leq C_{\divergence}\norm{\mathbf{u}}_{L_{\divergence}^2(\Omega)},
\end{equation}
with a constant $C_{\divergence}$ depending only on $\Omega$ (see e.g. \cite[Chap. III.2]{Galdi}).
Furthermore, we will use the following generalised Gagliardo-Nirenberg inequality:
\begin{lemma}\label{LEM_GAGNIR}
	Let $\Omega\subset \R^d,~d\in\N,$ be a bounded domain with Lipschitz-boundary and $f\in W^{m,r}\cap L^q,~1\leq q,r\leq \infty$. For any integer $j,~0\leq j<m$, suppose there is $\alpha\in\R$ such that 
	\begin{equation*}
	j - \frac{d}{p} = \left(m-\frac{d}{r}\right)\alpha + (1-\alpha)\left(-\frac{d}{q}\right),\quad \frac{j}{m}\leq \alpha\leq 1.
	\end{equation*}
	Then, there exists a positive constant $C$ depending only on $\Omega,d,m,j,q,r,$ and $\alpha$ such that
	\begin{equation}
	\label{LEM_GAGNIR_EST}\norml{p}{D^jf}\leq C\normw{m}{r}{f}^{\alpha}\norml{q}{f}^{1-\alpha}.
	\end{equation}
\end{lemma}
We will also need the following theorem concerning solvability of the divergence equation:
\begin{lemma}(\cite[Sec. III.3]{Galdi})\label{LEM_DIVEQU}
	Let $\Omega\subset\R^d,\,d\geq 2,$ be a bounded domain with Lipschitz-boundary and let $1<q<\infty$. Then, for every $f\in L^q$ and $\a\in \W^{1-1/q,q}(\del\Omega)$ satisfying
	\begin{equation}
	\label{DIV_COMP_COND}\intO f\dx  = \int_{\del\Omega}\a\cdot\n \d \mathcal{H}^{d-1},
	\end{equation}
	there exists at least one solution $\u \in \W^{1,q}$ of the problem
	\begin{alignat*}{3}
	\divergence(\u) &= f&&\quad\text{in }\Omega,\\
	\u  &= \a &&\quad \text{on }\del\Omega.
	\end{alignat*}
	In addition, the following estimate holds
	\begin{equation}
	\label{DIV_EST}\normW{1}{q}{\u}\leq C\left(\normL{q}{f}+\norm{\a}_{\W^{1-1/q,q}(\del\Omega)}\right),
	\end{equation}
	with $C$ depending only on $\Omega$ and $q$.
\end{lemma}

Finally, in the Galerkin ansatz (see Sec. 3) we will make use of the following lemma (see \cite{AbelsTerasawa} for a proof):
\begin{lemma}\label{LEM_STOKES}
	Let $\Omega\subset\R^d, d=2,3,$ be a bounded domain with $C^2$-boundary and outer unit normal $\n $ and $1<q<\infty$. Furthermore, assume that $g\in W^{1,q}$, $\f\in \L^q$, $c\in W^{1,r}$ with $r>d$, and the functions $\eta(\cdot),\,\lambda(\cdot)$ fulfil (A3) (see Assumptions \ref{ASS} below). Then, there exists a unique solution $(\v,p)\in \W^{2,q}\times W^{1,q}$ of the system 
	\begin{subequations}
		\label{STOKES_EQU}
		\begin{alignat}{3}\
		\label{STOKES_EQU_1}-\divergence(2\eta(c) \D\v +\lambda(c)\divergence(\v)\I )+\nu\v + \grad p &= \f&&\quad \text{a.\,e. in }\Omega,\\
		\label{STOKES_EQU_2}\divergence(\v) &= g&&\quad\text{a.\,e. in }\Omega,\\
		\label{STOKES_EQU_3}(2\eta(c) \D\v +\lambda(c)\divergence(\v)\I -p\I )\n  &= \mathbf{0}&&\quad \text{a.\,e. on }\del \Omega,
		\end{alignat}
	\end{subequations}
	satisfying the following estimate
	\begin{equation}
	\label{STOKES_EST}\normW{2}{q}{\v} + \normw{1}{q}{p}\leq C\left(\normL{q}{\f}+\normw{1}{q}{g}\right),
	\end{equation}
	with a constant $C$ depending only on $\eta_0,\eta_1,\lambda_0,q, \normw{1}{r}{c}$ and $\Omega$.
\end{lemma}

\section{Main results}
We make the following assumptions:
\begin{annahme} \label{ASS}\,
	\begin{enumerate}
		\item[(A1)] The positive constants $\epsilon,\,\nu,\,K,\,T$ are fixed,  $\chi$ is a fixed, non-negative constant. Furthermore, the function $\sigma_{\infty}\in L^{2}(L^2(\del\Omega))$ and the initial datum $\varphi_0\in H^1$ are prescribed.
		\item[(A2)] The mobility $m(\cdot)$ is continuous on $\R$ and satisfies
		\begin{equation*}
		m_0\leq m(t)\leq m_1\quad\forall t\in\R,
		\end{equation*}
		for positive constants $m_0,m_1$.
		\item[(A3)] The viscosities fulfil $\eta,\lambda\in C^2(\R)$ with bounded first derivatives and
		\begin{equation}
		\label{ASS_VISC}\eta_0\leq \eta(t)\leq \eta_1,\quad 0\leq \lambda(t)\leq \lambda_0\quad\forall t\in\R,
		\end{equation}
		for positive constants $\eta_0,\eta_1$ and a non-negative constant $\lambda_0$.
		\item[(A4)] The source terms are of the form
		\begin{equation*}
		\Gamma_{\v}(\varphi,\sigma) = b_{\v}(\varphi)\sigma + f_{\v}(\varphi),\quad \Gamma_{\varphi}(\varphi,\sigma) = b_{\varphi}(\varphi)\sigma + f_{\varphi}(\varphi),
		\end{equation*}
		where $b_{\v},f_{\v}\in C^1(\R)$ are bounded with bounded first derivatives and $b_{\varphi},f_{\varphi}\in C^0(\R)$ are bounded functions. The function $h\in C^0(\R)$ is continuous, bounded and non-negative.
		\item[(A5)] The function $\psi\in C^2(\R)$ is non-negative and can be written as
		\begin{equation}
		\label{ASS_PSI_1}\psi(s) = \psi_1(s) + \psi_2(s)\quad\forall s\in \R,
		\end{equation}
		where $\psi_1, \psi_2\in C^2(\R)$ and 
		\begin{alignat}{3}
		\label{ASS_PSI_2}R_1(1+|s|^{\rho-2})&\leq \psi_1''(s)\leq R_2(1+|s|^{\rho-2})&&\quad\forall s\in\R,\\
		\label{ASS_PSI_3}|\psi_2''(s)|&\leq R_3&&\quad\forall s\in\R,
		\end{alignat} 
		where $R_i,\,i=1,2,3,$ are positive constants with $R_1<R_2$ and $\rho\in [2,6]$.
		Furthermore, if $\rho=2$, we assume $2R_1>R_3$.
	\end{enumerate}
\end{annahme}
\begin{bemerkung}
	Using \textnormal{(A5)}, it is straightforward to check that there exist positive constants $R_i,\,i=4,5,6,$ such that 
	\begin{equation}
	\label{ASS_PSI_4}\psi(s)\geq R_4|s|^{\rho} - R_5\quad \forall s\in \R
	\end{equation}
	and
	\begin{equation}
	\label{ASS_PSI_5}|\psi'(s)|\leq R_6(1+|s|^{\rho-1})\quad \forall s\in \R.
	\end{equation}
	
\end{bemerkung}

We now introduce the weak formulation of \eqref{MEQ}, \eqref{BIC}:
\begin{definition}(Weak solution for \eqref{MEQ},\eqref{BIC}) \label{DEF_WSOL_1}We call a quintuple $(\varphi,\sigma,\mu,\v,p)$ a weak solution of \eqref{MEQ} and \eqref{BIC} if 
	\begin{align*}
	&\varphi \in H^1(0,T;(H^1)^*)\cap L^{\infty}(0,T;H^1)\cap L^2(0,T;H^3) ,\quad \mu\in L^2(0,T;H^1),\\
	&\sigma\in L^2(0,T;H^1),\quad \v\in L^2(0,T;\H^1),\quad p\in L^{2}(0,T;L^2),
	\end{align*}
	such that 
	\begin{equation*}
	\divergence(\v) = \Gamma_{\v}(\varphi,\sigma)\text{ a.\,e. in }Q,\quad \varphi(0)=\varphi_0\text{ a.\,e. in }\Omega,
	\end{equation*}
	and 
	\begin{subequations}
		\label{WFORM_1}
		\begin{align}
		\label{WFORM_1a}0 &= \intO \T(\v,p)\colon \grad\boldsymbol{\Phi}+\nu\v\cdot\boldsymbol{\Phi}  -  (\mu+\chi\sigma)\grad\varphi\cdot\boldsymbol{\Phi}\dx ,\\
		\label{WFORM_1b} 0 &= \langle\delt\varphi{,}\Phi\rangle_{H^1,(H^1)^*}  + \intO m(\varphi)\grad\mu\cdot\grad\Phi+ (\grad\varphi\cdot\v+ \varphi\Gamma_{\v}(\varphi,\sigma) - \Gamma_{\varphi}(\varphi,\sigma))\Phi\dx,  \\
		\label{WFORM_1c} 0 &= \intO (\mu+ \chi\sigma)\Phi -\epsilon^{-1}\Psi'(\varphi)\Phi - \epsilon\grad\varphi\cdot\grad\Phi \dx ,\\
		\label{WFORM_1d} 0  &=  \intO \grad\sigma\cdot\grad\Phi + h(\varphi)\sigma\Phi\dx +\int_{\del\Omega}K(\sigma-\sigma_{\infty})\Phi\d\mathcal{H}^{d-1} ,
		\end{align}
	\end{subequations}
	for a.\,e. $t\in(0,T)$ and for all $\boldsymbol{\Phi}\in \H^1,\,\Phi\in H^1$.
\end{definition}

\begin{theorem}\label{THM_WSOL_1}
	Let $\Omega\subset\R^d,\,d=2,3,$ be a bounded domain with $C^3$-boundary and assume that Assumptions \ref{ASS} is fulfilled. Then there exists a solution quintuple $(\varphi,\mu,\sigma,\v,p)$ of \eqref{MEQ}, \eqref{BIC} in the sense of Definition \ref{DEF_WSOL_1}. Furthermore, we have
	\begin{align*}
	&\varphi \in H^1(0,T;(H^1)^*)\cap L^{\infty}(0,T;H^1)\cap L^4(0,T;H^2)\cap L^2(0,T;H^3) ,\quad \sigma\in L^2(0,T;H^1),\\
	&\mu\in L^2(0,T;H^1),\quad\v\in L^2(0,T;\H^1),\quad p\in L^{2}(0,T;L^2),
	\end{align*}
	and the following estimate holds:
	\begin{align}
	\nonumber&\norm{\varphi}_{H^1((H^1)^*)\cap L^{\infty}(H^1)\cap L^4(H^2)\cap L^2(H^3) } + \norm{\sigma}_{L^2(H^1)} + \norm{\mu}_{L^2(H^1)}\\
	\label{THM_WSOL_1_EST_1}&\quad + \norm{\divergence(\varphi\v)}_{L^2(L^{\frac{3}{2}})}  + \norm{\v}_{L^2(\H^1)}+ \norm{p}_{L^2(L^2)}\leq C,
	\end{align}
	with a constant $C$ depending only on the system parameters, $\Omega$ and $T$.\newline
	If in addition $\sigma_{\infty}\in L^4(L^2(\del\Omega))$, we have
	\begin{equation*}
	\sigma\in L^4(0,T;H^1),\quad \mu\in L^4(0,T;L^2),\quad \v\in L^{\frac{8}{3}}(0,T;\H^1),
	\end{equation*}
	and
	\begin{equation}
	\label{THM_WSOL_1_EST_2} \norm{\sigma}_{L^4(H^1)} + \norm{\mu}_{L^4(L^2)}+ \norm{\v}_{L^{\frac{8}{3}}(\H^1)}
	+ \norm{\divergence(\varphi\v)}_{L^2(L^2)}  \leq C.
	\end{equation}
\end{theorem}
\begin{theorem}(The limit $K\to\infty$)\label{THM_KLIM} Let the assumptions of Theorem \ref{THM_WSOL_1} be fulfilled and assume in addition that $\sigma_{\infty}\in L^2(H^{\frac{1}{2}}(\del\Omega))$. Let $K>0$ and denote by $(\varphi_K,\mu_K,\sigma_K,\v_K,p_K)$ a weak solution of \eqref{MEQ}, \eqref{BIC} corresponding to $\varphi_0$ and $K$ in the sense of Definition \eqref{DEF_WSOL_1}. Then, as $K\to\infty$, we have
	\begin{alignat*}{3}
	\varphi_K&\to\varphi&&\quad\text{weakly-}*&&\quad\text{ in }H^1((H^1)^*)\cap L^{\infty}(H^1)\cap L^4(H^2)\cap L^2(H^3) ,\\
	\sigma_K&\to \sigma &&\quad\text{weakly-}*&&\quad\text{ in }L^2(H^1),\\
	\mu_K&\to \mu&&\quad\text{weakly}&&\quad\text{ in } L^2(H^1),\\
	p_K&\to p&&\quad\text{weakly}&&\quad\text{ in }L^2(L^2),\\
	\v_K&\to\v&&\quad\text{weakly}&&\quad\text{ in }L^2(\H^1),\\
	\divergence(\varphi_K\v_K)&\to \divergence(\varphi\v)&&\quad\text{weakly}&&\quad \text{ in }L^2(L^{\frac{3}{2}}),\\
	\sigma_K&\to \sigma_{\infty}&&\quad\text{strongly}&&\quad\text{ in }L^2(L^2(\del\Omega)),
	\end{alignat*}
	where $(\varphi,\mu,\sigma,\v,p)$ satisfies 
	\begin{equation*}
	\divergence(\v) = \Gamma_{\v}(\varphi,\sigma)\text{ a.\,e. in }Q,\quad \varphi(0)=\varphi_0\text{ a.\,e. in }\Omega,\quad \sigma\in (\sigma_{\infty}+L^2(0,T;H_0^1)),
	\end{equation*}
	and \eqref{WFORM_1} with \eqref{WFORM_1d} replaced by
	\begin{equation}
	\label{WF_KLIM_4}0 = \intO \grad\sigma\cdot\grad\xi + h(\varphi)\sigma\xi\dx ,
	\end{equation}
	for a.\,e. $t\in (0,T)$ and for all $\xi\in H_0^1$.
\end{theorem}
We now introduce the definition of weak solutions of the Cahn-Hilliard-Darcy system endowed with \eqref{BC_1}-\eqref{BC_2} and
\begin{equation}
\label{BC_3_DARCY}p = 0\quad\text{on } \Sigma.
\end{equation}
\begin{definition}\label{DEF_WS_DARCY} We call a quintuple $(\varphi,\mu,\sigma,\v,p)$ weak solution of the Cahn-Hilliard-Darcy system endowed with \eqref{BC_1}-\eqref{BC_2} and \eqref{BC_3_DARCY} if
	\begin{align*}
	&\varphi\in W^{1,\frac{8}{5}}(0,T;(H^1)^*)\cap L^{\infty}(0,T;H^1)\cap L^2(0,T;H^3),\quad\mu\in L^2(0,T;H^1),\\
	& \quad \sigma\in L^2(0,T;H^1),\quad \v\in L^2(0,T;\L^2),\quad p\in L^2(0,T;H_0^1),
	\end{align*}
	such that
	\begin{equation*}
	\varphi(0)=\varphi_0 \quad\text{a.\,e. in } \Omega
	\end{equation*}
	and
	\begin{subequations}\label{WF_DARCY}
		\begin{align}
		\label{WF_DARCY_a}0&= \langle \delt\varphi{,}\phi\rangle_{H^1} + \intO m(\varphi)\nabla\mu\cdot\nabla\phi\d x + \intO (\nabla\varphi\cdot\v + \varphi\Gamma_{\v}(\varphi,\sigma) - \Gamma_{\varphi}(\varphi,\sigma))\phi\d x,\\
		\label{WF_DARCY_b}0&= \intO(\mu+\chi\sigma)\phi - \epsilon^{-1}\psi'(\varphi)\phi - \epsilon\nabla\varphi\cdot\nabla\phi\d x,\\
		\label{WF_DARCY_c}0  &= \intO \nabla\sigma\cdot\nabla\phi + h(\varphi)\sigma\phi\d x + \int_{\del\Omega}K(\sigma-\sigma_{\infty})\phi\d\mathcal{H}^{d-1},\\
		\label{WF_DARCY_d}0 &= \intO \left(\v + \frac{1}{\nu}\nabla p-\frac{1}{\nu}(\mu+\chi\sigma)\nabla\varphi\right)\cdot\boldsymbol{\Phi}\d x,\\
		\label{WF_DARCY_e}0&= \intO \frac{1}{\nu}\nabla p\cdot\nabla\xi - \Gamma_{\v}(\varphi,\sigma)\xi - \frac{1}{\nu}(\mu+\chi\sigma)\nabla\varphi\cdot\nabla\xi\d x,
		\end{align}
	\end{subequations}
	for a.\,e. $t\in (0,T)$ and all $\phi\in H^1,\,\xi\in H_0^1,\,\boldsymbol{\Phi}\in\mathbf{L}^2$.
\end{definition}
The following theorem states that solutions of the Cahn-Hilliard-Darcy system can be found as the limit of the Cahn-Hilliard-Brinkman system when the viscosities tend to zero.
\begin{theorem}\label{THM_ZERO_VISC}
	Let $\Omega\subset\R^d,\,d=2,3,$ be a bounded domain with $C^3$-boundary and assume that \textnormal{(A1)}-\textnormal{(A2)}, \textnormal{(A4)}-\textnormal{(A5)} holds. Furthermore, let $\{\eta_n,\lambda_n\}_{n\in\N}$ be a sequence of function pairs fulfilling (A3) such that
	\begin{equation*}
	\norm{\eta_n(\cdot)}_{C^0(\R)}\to 0,\quad \norm{\lambda_n(\cdot)}_{C^0(\R)}\to 0\quad\text{as }n\to\infty.
	\end{equation*} 
	Let $(\varphi_n,\mu_n,\sigma_n,\v_n,p_n)$ be a sequence of weak solutions of the Cahn-Hilliard-Brinkman system in the sense of Definition \ref{DEF_WSOL_1} for $\eta(\cdot)=\eta_n(\cdot),\,\lambda(\cdot)= \lambda_n(\cdot)$ and originating from $\varphi_0\in H^1$. Then, at least for a subsequence, $(\varphi_n,\mu_n,\sigma_n,\v_n,p_n)$ converges to a weak solution $(\varphi,\mu,\sigma,\v,p)$ of the Cahn-Hilliard-Darcy system in the sense of Definition \ref{DEF_WS_DARCY} such that
	\begin{alignat*}{3}
	\varphi_{n}&\to\varphi&&\quad\text{weakly-}*&&\quad\text{ in } W^{1,\frac{8}{5}}((H^1)^*)\cap L^{\infty}(H^1)\cap L^4(H^2)\cap L^2(H^3),\\
	\sigma_{n}&\to \sigma &&\quad\text{weakly-}&&\quad\text{ in }L^2(H^1),\\
	\mu_{n}&\to \mu&&\quad\text{weakly}&&\quad\text{ in }L^2(H^1),\\
	p_{n}&\to p&&\quad\text{weakly}&&\quad\text{ in }L^{2}(L^2),\\
	\v_{n}&\to\v&&\quad\text{weakly}&&\quad\text{ in }L^2(\L^2)\cap L^2\left(L_{\divergence}^2(\Omega)\right),\\
	2\eta_n(\varphi_n) \D\v_{n}&\to \mathbf{0}&&\quad\text{weakly}&&\quad\text{ in }L^2(\L^2),\\
	\lambda_n(\varphi_n) \divergence(\v_n)\I&\to \mathbf{0}&&\quad\text{weakly}&&\quad\text{ in }L^2(\L^2),
	\end{alignat*}
	and
	\begin{equation*}
	\varphi_{n}\to\varphi\quad\text{strongly}\quad\text{ in }C^0(L^r)\cap L^2(W^{2,r})\text{ and a.\,e. in }Q,
	\end{equation*}
	for all $r\in [1,6)$. Moreover, it holds that
	\begin{alignat*}{3}
	\deln\varphi  &= p= 0 &&\quad \text{a.\,e. on }\Sigma,\\
	\divergence(\v) &= \Gamma_{\v}(\varphi,\sigma)&&\quad\text{a.\,e. in }Q,\\
	\v  &= \frac{1}{\nu}\left(-\nabla p + (\mu+\chi\sigma)\nabla\varphi\right)&&\quad\text{a.\,e. in }Q,\\
	\mu &= \epsilon^{-1}\psi'(\varphi)-\epsilon\laplace \varphi-\chi\sigma &&\quad\text{a.\,e. in }Q,
	\end{alignat*} 
	and
	\begin{align}
	\nonumber &\norm{\varphi}_{W^{1,\frac{8}{5}}((H^1)^*)\cap L^{\infty}(H^1)\cap L^4(H^2)\cap  L^2(H^3) } + \norm{\mu}_{L^2(H^1)} + \norm{\sigma}_{L^2(H^1)} \\
	\label{EST_DARCY}&\quad + \norm{\divergence(\varphi\v)}_{L^{\frac{8}{5}}((H^1)^*)}  + \norm{\v}_{L^2(L_{\divergence}^2(\Omega))}+ \norm{p}_{L^2(H_0^1)}\leq C,
	\end{align}	
	with a constant $C$ depending only on the system parameters and on $\Omega,\,T$.
\end{theorem}
To prove continuous dependence on the initial and boundary data, we make the following additional assumptions:
\begin{annahme}\label{ASS_1}\,
	\begin{enumerate}
		\item[(B1)] The mobility $m(\cdot)$ is a constant, without loss of generality we assume $m(\cdot)\equiv 1$.
		\item[(B2)] The functions $b_{\varphi}(\cdot)$,\,$f_{\varphi}(\cdot)$ and $h(\cdot),$ are Lipschitz continuous with Lipschitz constants $L_b$,\,$L_f$ and $L_h$, respectively.
		\item[(B3)] For $\psi'$ and $\psi''$, we assume that
		\begin{alignat}{3}
		\label{ASS_PSI_6}|\psi'(s_1)-\psi'(s_2)|&\leq k_1(1+|s_1|^4 + |s_2|^4)|s_1-s_2|&&\quad \forall s_1,s_2\in\R,\\
		\label{ASS_PSI_7}|\psi''(s_1)-\psi''(s_2)|&\leq k_2(1+|s_1|^3 + |s_2|^3)|s_1-s_2|&&\quad \forall s_1,s_2\in\R,
		\end{alignat}
		for some positive constants $k_1,k_2$.
	\end{enumerate}
\end{annahme}
\begin{theorem} \label{THM_CONTDEP}	Let $\Omega\subset\R^d,\,d=2,3,$ be a bounded domain with $C^3$-boundary and assume that Assumptions \ref{ASS} and \ref{ASS_1} hold. Then, for any two weak solution quintuples $\{\varphi_i,\mu_i,\sigma_i,\v_i,p_i\}$, $i=1,2,$ of \eqref{MEQ}, \eqref{BIC} satisfying 
	\begin{align*}
	&\varphi \in H^1(0,T;(H^1)^*)\cap L^{\infty}(0,T;H^1)\cap L^4(0,T;H^2)\cap L^2(0,T;H^3) ,\quad \sigma\in L^{4}(0,T;H^1),\\
	&\mu\in L^4(0,T;L^2)\cap L^2(0,T;H^1) ,\quad\v\in L^{\frac{8}{3}}(0,T;\H^1),\quad p\in L^{2}(0,T;L^2),
	\end{align*}
	with $\sigma_{i,\infty}\in L^4(L^2(\del\Omega))$ and $\varphi_i(0)=\varphi_{i,0}\in H^1$ for $i=1,2$, there exists a positive constant $C$ depending only on the system parameters and on $\Omega,T,L_h,L_b,L_f$ such that
	\begin{align}
	\nonumber &\sup_{t\in (0,T]}\left(\normh{1}{\varphi_1(t)-\varphi_2(t)}^2\right) + \norm{\varphi_1-\varphi_2}_{H^1((H^1)^*)\cap L^2(H^3)}^2 + \norm{\mu_1-\mu_2}_{L^2(H^1)}^2\\
	\nonumber&\quad + \norm{\sigma_1-\sigma_2}_{L^2(H^1)}^2 + \norm{\v_1-\v_2}_{L^2(\H^1)}^2 + \norm{p_1-p_2}_{L^2(L^2)}^2\\
	\label{THM_CONTDEP_EST}&\quad\leq C\left(\normh{1}{\varphi_{1,0}-\varphi_{2,0}}^2 + \norm{\sigma_{1,\infty}-\sigma_{2,\infty}}_{L^4(L^2(\del\Omega))}^2\right).
	\end{align}
\end{theorem}
We have the following notion of strong solutions:
\begin{definition}(Strong solution for \eqref{MEQ},\eqref{BIC}) \label{DEF_SSOL_1}We call a quintuple $(\varphi,\sigma,\mu,\v,p)$ a strong solution of \eqref{MEQ} and \eqref{BIC} if 
	\begin{align*}
	&\varphi \in H^1(0,T;L^2)\cap L^{\infty}(0,T;H^2)\cap L^2(0,T;H^4) ,\quad \mu\in L^2(0,T;H^2),\\
	&\sigma\in L^2(0,T;H^2),\quad\v\in L^2(0,T;\H^2),\quad p\in L^{2}(0,T;H^1),
	\end{align*}
	and \eqref{MEQ}, \eqref{BIC} are fulfilled almost everywhere in the respective sets.
\end{definition}
For the existence of strong solutions, we make the following additional assumptions:
\begin{annahme}\label{ASS_2}\,
	\item[(C1)] The mobility $m(\cdot)$ is a constant, without loss of generality we assume $m(\cdot)\equiv 1$. The function $h(\cdot)$ is Lipschitz continuous with Lipschitz constant $L_h$.
	\item[(C2)] The boundary datum $\sigma_{\infty}\in H^1(0,T;H^{\frac{1}{2}}(\del\Omega))$ and the initial datum $\varphi_0\in H_N^2$ are prescribed.
	\item[(C3)] The function $\psi\in C^3(\R)$ fulfils 
	\begin{equation}
	\label{ASS_PSI_8}|\psi'''(s)|\leq k_3(1+|s|^3)\quad\forall s\in \R,
	\end{equation}
	for a positive constant $k_3$.
\end{annahme}
We have the following result concerning strong solutions:
\begin{theorem}\label{THM_SSOL}
	Let $\Omega\subset\R^d,\,d=2,3,$ be a bounded domain with $C^4$-boundary and assume that Assumptions \ref{ASS} and \ref{ASS_2} hold. Then there exists a solution quintuple $(\varphi,\mu,\sigma,\v,p)$ of \eqref{MEQ}, \eqref{BIC} in the sense of Definition \ref{DEF_SSOL_1}. Furthermore, we have
	\begin{equation*}
	\varphi \in C^0(\bar{Q}),\quad \mu\in L^{\infty}(0,T;L^2),\quad \sigma\in  H^1(0,T;H^1)\cap L^{\infty}(0,T;H^2),\quad\v\in L^{8}(0,T;\H^2),\quad p\in L^{8}(0,T;H^1),
	\end{equation*}
	and the following estimate holds:
	\begin{align}
	\nonumber&\norm{\varphi}_{H^1(L^2)\cap C^0(\bar{Q})\cap L^{\infty}(H^2)\cap L^2(H^4) } + \norm{\sigma}_{H^1(H^1)\cap L^{\infty}(H^2)} + \norm{\mu}_{L^{\infty}(L^2)\cap L^2(H^2)}\\
	\label{THM_SSOL_EST}&\quad + \norm{\divergence(\varphi\v)}_{L^2(L^2)}  + \norm{\v}_{L^{8}(\H^2)}+ \norm{p}_{L^8(H^1)}\leq C.
	\end{align}
\end{theorem}

\section{Existence of weak solutions} 
In order to prove the result, we now derive a priori estimates for \eqref{WFORM_1a}-\eqref{WFORM_1d}. By $C$, we denote a positive constant not depending on $(\varphi,\mu,\sigma,\v,p)$ which may vary from line to line. The duality pairing in \eqref{WFORM_1b} can be replaced by the $L^2$-product for smooth enough functions which is satisfied for example by a Galerkin-ansatz. Approximating solutions can be constructed by applying a Galerkin approximation with respect to $\varphi$ and $\mu$ and at the same time solving for $\v,\,p$ and $\sigma$ in the corresponding whole function spaces (for details, we refer to \cite{EbenbeckGarcke}, \cite{GarckeLam1}). In the following, we will write $\Gamma_{\varphi},\,\Gamma_{\v},$ instead of $\Gamma_{\varphi}(\varphi,\sigma),\,\Gamma_{\v}(\varphi,\sigma)$.
\subsection{A-priori-estimates}
\subsubsection{Estimating the nutrient concentration}
Choosing $\Phi = \sigma$ in \eqref{WFORM_1d} and using the non-negativity of $h(\cdot)$, we obtain
\begin{equation}
\label{APRI_EQ_1}\intO |\grad\sigma|^2\dx   + K\int_{\del\Omega}|\sigma|^2 \d\mathcal{H}^{d-1} \leq  K\int_{\del\Omega}\sigma\sigma_{\infty}\d\mathcal{H}^{d-1}. 
\end{equation}
Using Hölder's and Young's inequalities, we have
\begin{equation}
\label{APRI_EQ_2}\left|K\int_{\del\Omega}\sigma\sigma_{\infty}\d\mathcal{H}^{d-1}\right|\leq \frac{K}{2}\norm{\sigma}_{L^2(\del\Omega)}^2 + \frac{K}{2}\norm{\sigma_{\infty}}_{L^2(\del\Omega)}^2.
\end{equation}
Using \eqref{APRI_EQ_1}-\eqref{APRI_EQ_2} and Poincaré's inequality, we deduce that
\begin{equation}
\label{APRI_EQ_3}\normh{1}{\sigma}\leq C\norm{\sigma_{\infty}}_{L^2(\del\Omega)}.
\end{equation}
Moreover, by the continuous embedding $H^1\hookrightarrow L^p,\, p\in[2,6],$ and \textnormal{(A4)}, we have
\begin{equation}
\label{APRI_EQ_4}\norml{p}{\Gamma_{\varphi}} + \norml{p}{\Gamma_{\v}}\leq C\left(1+\norm{\sigma_{\infty}}_{L^2(\del\Omega)}\right)\quad\forall p\in[2,6].
\end{equation}
\subsubsection{An energy identity}
Due to \textnormal{(A4)}, \eqref{APRI_EQ_3}-\eqref{APRI_EQ_4}, there exists a solution $\u \in \H^1$  of the problem
\begin{equation}
\label{APRI_EQ_6}\divergence(\u) = \Gamma_{\v} \quad\text{in }\Omega,\qquad \u  = \frac{1}{|\del\Omega|}\left(\intO\Gamma_{\v}\dx \right)\n \eqqcolon \a\quad\text{on }\del\Omega,
\end{equation}
satisfying the estimate
\begin{equation}
\label{APRI_EQ_7}\norm{\u }_{\H^1}\leq C\norml{2}{\Gamma_{\v}}\leq C\left(1+\norm{\sigma_{\infty}}_{L^2(\del\Omega)}\right).
\end{equation}
Choosing $\boldsymbol{\Phi}=\v-\u$ in \eqref{WFORM_1a}, $\Phi = \mu+\chi\sigma$ in \eqref{WFORM_1b}, $\Phi = \delt\varphi$ in \eqref{WFORM_1c},  and summing the resulting identities, we obtain
%
\begin{align}
\nonumber &\frac{\d}{\dt } \intO \epsilon^{-1}\psi(\varphi)+\frac{\epsilon}{2}|\grad\varphi|^2 \dx  +\intO m(\varphi)|\grad\mu|^2 +\intO 2\eta(\varphi)|\D\v|^2 +\nu|\v|^2\dx \\
\nonumber&\quad = \intO -m(\varphi)\chi\grad\mu\cdot\grad\sigma + (\Gamma_{\varphi}-\varphi\Gamma_{\v})(\mu+\chi\sigma)\dx \\
\label{APRI_EQ_8}& \quad +\intO 2\eta(\varphi) \D\v\colon \grad\u  + \nu\v\cdot \u \dx  - \intO (\mu + \chi\sigma)\grad\varphi \cdot \u \dx.
\end{align}

\subsubsection{Estimating the right hand side of the energy identity}
Using Hölder's and Young's inequalities together with \eqref{APRI_EQ_7}, we have
\begin{equation}
\label{APRI_EQ_9}\left|\intO 2\eta(\varphi) \D\v \colon\grad\u  + \nu\v \cdot\u \dx \right| \leq \normL{2}{\sqrt{\eta(\varphi)}\D\v}^2 + \frac{\nu}{2}\normL{2}{\v}^2 + C\left(1+\norml{\infty}{\eta(\varphi)}\right)\norml{2}{\Gamma_{\v}}^2.
\end{equation}
We now want to estimate the terms involving $\Gamma_{\v}$ and $\Gamma_{\varphi}$.
Using Hölder's, Poincaré's and Young's inequalities together with \eqref{APRI_EQ_3}-\eqref{APRI_EQ_4}, we obtain
\begin{align}
\nonumber \left|\intO  \Gamma_{\varphi}(\mu+\chi\sigma)\dx \right| &\leq C_P\norml{2}{\Gamma_{\varphi}}\big(|(\mu+\chi\sigma)_{\Omega}| + \normL{2}{\grad(\mu+\chi\sigma)}\big)\\
\label{APRI_EQ_10}&\leq C\left(1+\norm{\sigma_{\infty}}_{L^2(\del\Omega)}^2\right)\big(1+|(\mu+\chi\sigma)_{\Omega}|\big) + \frac{m_0}{8}\normL{2}{\grad\mu}^2.
\end{align}
With similar arguments, using the continuous embedding $H^1\hookrightarrow L^6$, we obtain
\begin{align}
\nonumber \left|\intO  \Gamma_{\v}\varphi(\mu+\chi\sigma)\dx \right| &\leq C\norml{3}{\Gamma_{\v}}\norml{2}{\varphi}\norml{6}{\mu+\chi\sigma} \\ 
\nonumber &\leq C\norml{3}{\Gamma_{\v}}\norml{2}{\varphi}\big(|(\mu+\chi\sigma)_{\Omega}| + \normL{2}{\grad(\mu+\chi\sigma)}\big)\\
\label{APRI_EQ_11} &\leq C\left(1+\norm{\sigma_{\infty}}_{L^2(\del\Omega)}^2\right)\big(1+\norml{2}{\varphi}^2 + \norml{2}{\varphi}|(\mu+\chi\sigma)_{\Omega}|\big) + \frac{m_0}{8}\normL{2}{\grad\mu}^2.
\end{align}
Now, choosing $\Phi = 1$ in \eqref{WFORM_1c} and using \eqref{ASS_PSI_5}, we obtain
\begin{equation*}
\left|\intO \mu+\chi\sigma\dx \right| = \left|\intO \epsilon^{-1}\psi'(\varphi)\dx\right|\leq \epsilon^{-1}R_6\intO 1+|\varphi|^{\rho-1}\dx \leq C\left(1+\norml{\rho-1}{\varphi}^{\rho-1}\right),
\end{equation*}
hence
\begin{equation*}
|(\mu+\chi\sigma)_{\Omega}|\leq C\left(1+\norml{\rho-1}{\varphi}^{\rho-1}\right)\leq C\left(1+\norml{\rho}{\varphi}^{\rho-1}\right).
\end{equation*}
In particular, using Young's inequality, the continuous embedding $L^{\rho}\hookrightarrow L^2,\,\rho\in [2,6],$ and \eqref{ASS_PSI_4}, this implies
\begin{align}
\label{APRI_EQ_12} |(\mu+\chi\sigma)_{\Omega}|&\leq C\left(1+\norml{\rho}{\varphi}^{\rho-1}\right) \leq C\left(1+\norml{\rho}{\varphi}^{\rho}\right)\leq C\left(1+\norml{1}{\psi(\varphi)}\right),\\
\label{APRI_EQ_13} \norml{2}{\varphi}|(\mu+\chi\sigma)_{\Omega}|&\leq C\left(\norml{2}{\varphi} + \norml{\rho}{\varphi}^{\rho}\right)\leq C\left(1+\norml{\rho}{\varphi}^{\rho}\right)\leq C\left(1+\norml{1}{\psi(\varphi)}\right).
\end{align}
Using \eqref{APRI_EQ_12}-\eqref{APRI_EQ_13} in \eqref{APRI_EQ_10}-\eqref{APRI_EQ_11} and applying the continuous embedding $L^{\rho}\hookrightarrow L^2,\,\rho\in [2,6],$ together with \eqref{ASS_PSI_4}, we end up with
\begin{equation}
\label{APRI_EQ_14} \left|\intO  (\Gamma_{\varphi}-\varphi\Gamma_{\v})(\mu+\chi\sigma) \dx\right| \leq C\left(1+\norm{\sigma_{\infty}}_{L^2(\del\Omega)}^2\right)\left(1+\norml{1}{\psi(\varphi)}\right) + \frac{m_0}{4}\normL{2}{\grad\mu}^2.
\end{equation}
For the first term on the r.h.s. of \eqref{APRI_EQ_8}, applying Hölder's and Young's inequalities, \textnormal{(A2)} and \eqref{APRI_EQ_3}, we obtain
\begin{equation}
\label{APRI_EQ_15} \left|\intO m(\varphi)\chi\grad\mu\cdot\grad\sigma\dx \right| \leq m_1\chi\normL{2}{\grad\mu}\normL{2}{\grad\sigma}\leq \frac{m_0}{8}\normL{2}{\grad\mu}^2 + C\norm{\sigma_{\infty}}_{L^2(\del\Omega)}^2.
\end{equation}
\subsubsection{Estimating the remaining term} 
\paragraph{The case $\rho=2$:}Using Hölder's, Young's and Poincaré's inequalities, the continuous embeddings $\H^1\hookrightarrow \L^6,\, H^1\hookrightarrow L^3$ and \eqref{ASS_PSI_4}, \eqref{APRI_EQ_7}, \eqref{APRI_EQ_12}, we obtain
\begin{align}
\nonumber\left|\intO (\mu+\chi\sigma)\grad\varphi\cdot\u\dx \right| &\leq \norml{3}{\mu+\chi\sigma}\normL{2}{\grad\varphi}\normL{6}{\u}\\
\nonumber&\leq C\big(|(\mu+\chi\sigma)_{\Omega}| + \normL{2}{\grad(\mu+\chi\sigma)}\big)\normL{2}{\grad\varphi}\normH{1}{\u}\\
\nonumber&\leq C\big(1+\norml{2}{\varphi} + \normL{2}{\grad(\mu+\chi\sigma)}\big)\normL{2}{\grad\varphi}\left(1+\norm{\sigma_{\infty}}_{L^2(\del\Omega)}\right)\\
\nonumber&\leq C\left(1+\norm{\sigma_{\infty}}_{L^2(\del\Omega)}^2\right)\left(1+\norml{2}{\varphi}^2 + \normL{2}{\grad\varphi}^2\right) + \frac{m_0}{8}\normL{2}{\grad\mu}^2\\
\label{APRI_EQ_16}&\leq  C\left(1+\norm{\sigma_{\infty}}_{L^2(\del\Omega)}^2\right)\left(1+\norml{1}{\psi(\varphi)} + \normL{2}{\grad\varphi}^2\right) + \frac{m_0}{8}\normL{2}{\grad\mu}^2.
\end{align}
\paragraph{The case $\rho\in (2,6]$:} In this case, we need a more subtle argument. Choosing $\Phi = -\laplace\varphi$ in \eqref{WFORM_1c}, integrating by parts and using \eqref{ASS_PSI_1}, it holds
\begin{equation}
\label{APRI_EQ_17}\epsilon\intO |\laplace\varphi|^2\dx + \epsilon^{-1}\intO \psi_1''(\varphi)|\grad\varphi|^2\dx = -\epsilon^{-1}\intO \psi_2''(\varphi)|\grad\varphi|^2\dx + \intO \grad(\mu+\chi\sigma)\cdot\grad\varphi\dx.
\end{equation}
Neglecting the non-negative term $\epsilon\intO |\laplace\varphi|^2\dx$ on the l.h.s. of this equation and using \eqref{ASS_PSI_2}-\eqref{ASS_PSI_3} together with Young's inequality, we obtain
\begin{equation}
\label{APRI_EQ_18}\intO |\varphi|^{\rho-2}|\grad\varphi|^2\dx \leq \left(1+\frac{R_3}{R_1} + \frac{\epsilon^2}{4\delta R_1^2}\right)\normL{2}{\grad\varphi}^2 + \delta\normL{2}{\grad(\mu+\chi\sigma)}^2,
\end{equation}
with $\delta>0$ to be chosen later. Observing that
\begin{equation*}
\left| \grad\left(\frac{2|\varphi|^{\frac{\rho}{2}}}{\rho}\right)\right| = |\varphi|^{\frac{\rho-2}{2}}|\grad\varphi|,
\end{equation*}
from \eqref{APRI_EQ_18} we deduce that
\begin{equation}
\label{APRI_EQ_19}\normL*{2}{\grad(|\varphi|^{\frac{\rho}{2}})}^2\leq \frac{\rho^2}{4}\left(1+\frac{R_3}{R_1} + \frac{\epsilon^2}{4\delta R_1^2}\right)\normL{2}{\grad\varphi}^2 + \frac{\delta\rho^2}{4}\normL{2}{\grad(\mu+\chi\sigma)}^2.
\end{equation}
Now, applying the trace theorem and \eqref{APRI_EQ_19} yields
\begin{align*}
\nonumber \norm*{|\varphi|^{\frac{\rho}{2}}}_{L^2(\del\Omega)}^2 &\leq C_{tr}^2\left(\norml*{2}{|\varphi|^{\frac{\rho}{2}}}^2 + \normL*{2}{\grad\left(|\varphi|^{\frac{\rho}{2}}\right)}^2\right)\\
&\leq C_{tr}^2\left(\norml{\rho}{\varphi}^{\rho} +\frac{\rho^2}{4}\left(1+\frac{R_3}{R_1} + \frac{\epsilon^2}{4\delta R_1^2}\right)\normL{2}{\grad\varphi}^2 + \frac{\delta\rho^2}{4}\normL{2}{\grad(\mu+\chi\sigma)}^2\right),
\end{align*}
hence
\begin{equation}
\label{APRI_EQ_20} \norm{\varphi}_{L^{\rho}(\del\Omega)}^{\rho}\leq C_{tr}^2\left(\norml{\rho}{\varphi}^{\rho} +\frac{\rho^2}{4}\left(1+\frac{R_3}{R_1} + \frac{\epsilon^2}{4\delta R_1^2}\right)\normL{2}{\grad\varphi}^2 + \frac{\delta\rho^2}{4}\normL{2}{\grad(\mu+\chi\sigma)}^2\right).
\end{equation}
Now, upon integrating by parts and recalling $\eqref{APRI_EQ_6}_2$,
we calculate
\begin{align}
\nonumber\intO (\mu+\chi\sigma)\grad\varphi\cdot\u\dx &= \intO \big(\mu+\chi\sigma - (\mu+\chi\sigma)_{\Omega}\big)\grad\varphi\cdot\u\dx + (\mu+\chi\sigma)_{\Omega}\intO \grad\varphi\cdot\u\dx\\
\nonumber&= \intO \big(\mu+\chi\sigma - (\mu+\chi\sigma)_{\Omega}\big)\grad\varphi\cdot\u\dx \\
\label{APRI_EQ_21} &\quad + (\mu+\chi\sigma)_{\Omega}\left(\frac{1}{|\del\Omega|}\intO \Gamma_{\v}\dx \int_{\del\Omega}\varphi\d\mathcal{H}^{d-1} - \intO \varphi\Gamma_{\v}\dx\right).
\end{align}
Using Hölder's, Young's and Poincaré's inequalities, the continuous embedding $\H^1\hookrightarrow\L^6$ and \eqref{APRI_EQ_7}, it is straightforward to check that
\begin{equation}
\label{APRI_EQ_22}\left|\intO \big(\mu+\chi\sigma - (\mu+\chi\sigma)_{\Omega}\big)\grad\varphi\cdot\u\dx\right| \leq C_{\delta_1}(1+\norm{\sigma_{\infty}}_{L^2(\del\Omega)}^2)(1+\normL{2}{\grad\varphi}^2) + \delta_1 \normL{2}{\grad\mu}^2,
\end{equation}
with $\delta_1>0$ to be chosen. Using \eqref{APRI_EQ_7}, \eqref{APRI_EQ_13} and Hölder's inequality, we obtain
\begin{equation}
\label{APRI_EQ_23} \left|(\mu+\chi\sigma)_{\Omega}\intO \varphi\Gamma_{\v}\dx\right| \leq |(\mu+\chi\sigma)_{\Omega}| \norml{2}{\Gamma_{\v}}\norml{2}{\varphi}\leq C\left(1+\norm{\sigma_{\infty}}_{L^2(\del\Omega)}\right)\left(1+\norml{1}{\psi(\varphi)}\right)
\end{equation}
Now, using Hölder's and Young's inequalities, \eqref{ASS_PSI_4}, \eqref{APRI_EQ_7}, \eqref{APRI_EQ_12}, \eqref{APRI_EQ_20} and recalling $\rho >2$, we deduce that
\begin{align}
\nonumber  \left|(\mu+\chi\sigma)_{\Omega}\frac{1}{|\del\Omega|}\intO \Gamma_{\v}\dx \int_{\del\Omega}\varphi\d\mathcal{H}^{d-1}\right| &\leq C|(\mu+\chi\sigma)_{\Omega}| \norml{2}{\Gamma_{\v}}\norm{\varphi}_{L^{\rho}(\del\Omega)}\\
\nonumber &\leq C\left(1+\norml{\rho}{\varphi}^{\rho-1}\right)\norml{2}{\Gamma_{\v}}\norm{\varphi}_{L^{\rho}(\del\Omega)}\\
\nonumber &\leq C_{\delta_2}(1+\norml{\rho}{\varphi}^{\rho})\norml{2}{\Gamma_{\v}}^{\frac{\rho}{\rho-1}} + \delta_2\norm{\varphi}_{L^{\rho}(\del\Omega)}^{\rho}\\
\label{APRI_EQ_24}&\leq C_{\delta_2}\left(1+\norm{\sigma_{\infty}}_{L^2(\del\Omega)}^2\right)(1+\norml{1}{\psi(\varphi)})+ \delta_2\norm{\varphi}_{L^{\rho}(\del\Omega)}^{\rho}
\end{align}
where we used that $\frac{\rho-1}{\rho} + \frac{1}{\rho} = 1$. Now, plugging in \eqref{APRI_EQ_22}-\eqref{APRI_EQ_24} into \eqref{APRI_EQ_21}, using \eqref{ASS_PSI_4}, \eqref{APRI_EQ_20} and choosing $\delta,\delta_1,\delta_2$ small enough, we finally obtain
\begin{equation}
\label{APRI_EQ_25} \left|\intO (\mu+\chi\sigma)\grad\varphi\cdot\u\dx\right|\leq \left(1+\norm{\sigma_{\infty}}_{L^2(\del\Omega)}^2\right)\left(1+\norml{1}{\psi(\varphi)} + \normL{2}{\grad\varphi}^2\right) + \frac{m_0}{8}\normL{2}{\grad\mu}^2.
\end{equation}
Plugging in \eqref{APRI_EQ_9}, \eqref{APRI_EQ_14}-\eqref{APRI_EQ_16}, \eqref{APRI_EQ_25}  into \eqref{APRI_EQ_8}, using \eqref{ASS_PSI_4} and \eqref{APRI_EQ_3}, we end up with
\begin{align*}
 &\frac{\d}{\dt } \intO \epsilon^{-1}\psi(\varphi)+\frac{\epsilon}{2}|\grad\varphi|^2 \dx  + \frac{m_0}{4}\normL{2}{\grad\mu}^2+  \frac{\nu}{2}\normL{2}{\v}^2 + \normL{2}{\sqrt{\eta(\varphi)}\D\v}^2 + \normh{1}{\sigma}^2\\
&\leq \alpha(t)\left(1+ \normL{2}{\grad\varphi(t)}^2 + \norml{1}{\psi(\varphi(t))}\right),
\end{align*}
where, recalling \textnormal{(A1)} and \textnormal{(A3)},
\begin{equation}
\label{APRI_EQ_25a}\alpha(t)\coloneqq C\left(1+\norml{\infty}{\eta(\varphi(t))}\right)\left(1 + \norm{\sigma_{\infty}(t)}_{L^2(\del\Omega)}^2\right)\in L^1(0,T).\\
\end{equation}
Integrating the last inequality in time from $0$ to $s\in (0,T]$
and applying Gronwall's Lemma (see \cite[Lemma 3.1]{GarckeLam3}) yields
\begin{align}
\nonumber &  \epsilon^{-1}\norml{1}{\psi(\varphi(s))}+\frac{\epsilon}{2}\normL{2}{\grad\varphi(s)}^2 + \int_{0}^{s}\frac{m_0}{4}\normL{2}{\grad\mu}^2+  \frac{\nu}{2}\normL{2}{\v}^2 + \normL{2}{\sqrt{\eta(\varphi)}\D\v}^2 + \normh{1}{\sigma}^2\d t\\
\label{APRI_EQ_26}&\leq \left(\epsilon^{-1}\norml{1}{\psi(\varphi_0)}+\frac{\epsilon}{2}\normL{2}{\grad\varphi_0}^2 + \int_{0}^{s}\alpha(t)\dt\right)\exp\left(\int_{0}^{s}\alpha(t)\dt\right)\quad\forall s\in (0,T].
\end{align}
Due to \textnormal{(A1)}, \textnormal{(A5)} and the continuous embedding $H^1\hookrightarrow L^6$, we have $\psi(\varphi_0)\in L^1,\,\varphi_0\in \L^2$. Then, due to Korn's inequality (see \cite[Thm. 6.3-3]{Ciarlet}) and \textnormal{(A3)}, taking the supremum over all $s\in (0,T]$ in \eqref{APRI_EQ_26} implies
\begin{equation}
\label{APRI_EQ_27}\esssup_{s\in (0,T]}\left(\norml{1}{\psi(\varphi(s))} + \normL{2}{\grad\varphi(s)}^2\right) +\int_{0}^{T}\normH{1}{\v}^2 + \normL{2}{\grad\mu}^2 + \normh{1}{\sigma}^2 \dt \leq C.
\end{equation}
Recalling \eqref{ASS_PSI_4}, using Poincaré's inequality, \eqref{APRI_EQ_12}, \eqref{APRI_EQ_27} and the fact that $\rho\geq 2$, this in particular gives
\begin{equation*}
\esssup_{s\in (0,T]}\normh{1}{\varphi(s)}^2 + \intT \normh{1}{\mu}^2\dt\leq C.
\end{equation*}
Combining the last inequality with \eqref{APRI_EQ_27} yields
\begin{equation}
\label{APRI_EQ_28}\esssup_{s\in (0,T]}\left(\norml{1}{\psi(\varphi(s))} + \normh{1}{\varphi(s)}^2\right) +\int_{0}^{T}\normH{1}{\v}^2 + \normh{1}{\mu}^2 + \normh{1}{\sigma}^2 \dt \leq C.
\end{equation}
Due to \textnormal{(A4)} and the continuous embedding $H^1\hookrightarrow L^6$, this implies
\begin{equation}
\label{APRI_EQ_29}\norm{\Gamma_{\v}}_{L^2(L^6)} + \norm{\Gamma_{\varphi}}_{L^2(L^6)}\leq C.
\end{equation}

\subsubsection{Estimating the pressure}
Using Lemma \ref{LEM_DIVEQU}, we deduce that there is at least one solution $\q\in \H^1$ of the system
\begin{equation*}
\divergence(\q) = p\quad\text{in }\Omega,\qquad
\q = \frac{1}{|\del\Omega|}\left(\intO p\dx \right)\n \quad\text{on }\del\Omega,
\end{equation*}
such that
\begin{equation}
\label{APRI_EQ_30}\normH{1}{\q}\leq C_d\norml{2}{p},
\end{equation}
with $C_d$ depending only on $\Omega$ and $q=2$. Notice that the compatibility condition \eqref{DIV_COMP_COND} is satisfied since
\begin{equation*}
\int_{\del\Omega} \q\cdot\n \d \mathcal{H}^{d-1} = \frac{1}{|\del\Omega|}\left(\intO p\dx \right)\int_{\del\Omega}\n \cdot\n \d \mathcal{H}^{d-1} = \intO p\dx .
\end{equation*}
Choosing $\boldsymbol{\Phi}=\q$ in \eqref{WFORM_1a} and using $\divergence(\v)=\Gamma_{\v}$ a.\,e. in $Q$, we obtain
\begin{equation}
\label{APRI_EQ_31}\intO  |p|^2\dx  = \intO (2\eta(\varphi) \D\v + \lambda(\varphi)\Gamma_{\v}\I )\colon \grad\q\dx + \intO (\nu\v - (\mu + \chi\sigma)\grad\varphi)\cdot\q\dx . 
\end{equation}
Using \eqref{APRI_EQ_30}, \textnormal{(A3)} and Hölder's and Young's inequalities, an easy calculation shows that
\begin{equation}
\label{APRI_EQ_32}\norml{2}{p}^2 \leq C\left(\norml{\infty}{\eta(\varphi)}\normL{2}{\sqrt{\eta(\varphi)}\D\v}^2 + \norml{\infty}{\lambda(\varphi)}^2\norml{2}{\Gamma_{\v}}^2 + \normL{2}{\v}^2 + \norml{3}{\mu + \chi\sigma}^2\normL{2}{\grad\varphi}^2\right).
\end{equation}
Integrating this inequality in time from $0$ to $T$ and using Hölder's inequality, we arrive at
\begin{align*}
\norm{p}_{L^2(L^2)}^2&\leq C\left(\norm{\eta(\cdot)}_{L^{\infty}(\R)}\norm{\sqrt{\eta(\varphi)}\D\v}_{L^2(\L^2)}^2 + \norm{\lambda(\cdot)}_{L^{\infty}(\R)}^2\norm{\Gamma_{\v}}_{L^2(L^2)}^2\right)\\
&\quad + C\left(\norm{\v}_{L^2(\L^2)}^2+ \norm{\mu+\chi\sigma}_{L^2(L^3)}^2\norm{\grad\varphi}_{L^{\infty}(\L^2)}^2\right).
\end{align*}
Using \eqref{APRI_EQ_26}, \eqref{APRI_EQ_28}-\eqref{APRI_EQ_29} and \textnormal{(A3)}, the last inequality implies
\begin{equation}
\label{APRI_EQ_33} \norm{p}_{L^2(L^2)}\leq C.
\end{equation}

\subsubsection{Higher order estimates for $\varphi$}
Our aim is to show that
\begin{equation}
\label{APRI_EQ_34}\norm{\varphi}_{L^2(H^3)}\leq C.
\end{equation}
Choosing $\Phi = -\laplace\varphi$ in \eqref{WFORM_1c}, integrating by parts and neglecting the non-negative term resulting from $\psi_1$ (see \textnormal{(A5)}), we obtain 
\begin{equation*}
\epsilon\norml{2}{\laplace\varphi}^2 + \intO \epsilon^{-1}\psi_2''(\varphi)|\grad\varphi|^2\dx \leq  \intO \grad(\mu + \chi\sigma)\cdot\grad\varphi \dx .
\end{equation*}
Using Hölder's inequality and the assumptions on $\psi_2$, we therefore get
\begin{equation*}
\label{higher_order_estimates_4}\epsilon\norml{2}{\laplace\varphi}^2 \leq \epsilon^{-1}R_3\normL{2}{\grad\varphi}^2 + \normL{2}{\grad(\mu + \chi\sigma)}\normL{2}{\grad\varphi}.
\end{equation*}
Taking the square of this inequality and integrating in time from $0$ to $T$, we obtain
\begin{align}
\nonumber \epsilon^2\int_{0}^{T}\norml{2}{\laplace\varphi}^4 &\leq C\intT \normL{2}{\grad\varphi}^4 + \normL{2}{\grad\varphi}^2\normL{2}{\grad(\mu + \chi\sigma)}^2\dt \\
\nonumber &\leq C\left(\norm{\grad\varphi}_{L^{\infty}(\L^2)}^4 + \norm{\grad(\mu +\chi\sigma)}_{L^2(\L^2)}^2\norm{\grad\varphi}_{L^{\infty}(\L^2)}^2\right).
\end{align}
Applying elliptic regularity theory and \eqref{APRI_EQ_28}, this gives
\begin{equation}
\label{APRI_EQ_35}\norm{\varphi}_{L^4(H^2)}\leq C.
\end{equation}
Next, we take $\Phi = \laplace^2\varphi$ in \eqref{WFORM_1c}, integrate by parts and in time from $0$ to $T$ to obtain
\begin{equation}
\label{APRI_EQ_36}\epsilon\norm{\grad\laplace\varphi}_{L^2(\L^2)}^2 =  - \intT\intO \grad(\mu + \chi\sigma)\cdot\grad\laplace\varphi\dx \dt + \intT\intO \epsilon^{-1}\psi''(\varphi)\grad\varphi\cdot\grad\laplace\varphi\dx \dt .
\end{equation}
For the first term on the r.h.s. of \eqref{APRI_EQ_36}, applying Hölder's and Young's inequality gives
\begin{equation}
\label{APRI_EQ_37}\left|\intT\intO \grad(\mu +  \chi\sigma)\cdot\grad\laplace\varphi\dx \dt \right| \leq C \left(\norm{\mu}_{L^2(H^1)}^2 + \norm{\sigma}_{L^2(H^1)}^2\right) + \frac{\epsilon}{4}\norm{\grad\laplace\varphi}_{L^2(\L^2)}^2.
\end{equation}
Due to \eqref{LEM_GAGNIR_EST}, \eqref{ASS_PSI_2}-\eqref{ASS_PSI_3} and \eqref{APRI_EQ_28}, using Hölder's and Young's inequalities yields
\begin{align}
\nonumber \left|\intT\intO \epsilon^{-1}\psi''(\varphi)\grad\varphi\cdot\grad\laplace\varphi\dx \dt\right| &\leq C \intT\intO \left(1+|\varphi|^{\rho-2}\right)|\nabla\varphi||\nabla\laplace\varphi|\d x\d t\\
\nonumber &\leq C\intO \left(1+\norml{\infty}{\varphi}^{\rho-2}\right)\normL{2}{\nabla\varphi}\normL{2}{\nabla\laplace\varphi}\d t\\
\nonumber &\leq C\intT \left(1+\normh{2}{\varphi}^{\frac{\rho-2}{2}}\right)\normL{2}{\nabla\laplace\varphi}\d t\\
\label{APRI_EQ_38}&\leq C\left(1+\norm{\varphi}_{L^{\rho-2}(H^2)}^{\rho-2}\right) + \frac{\epsilon}{4}\norm{\grad\laplace\varphi}_{L^2(\L^2)}^2.
\end{align}
Recalling that $\rho-2\leq 4$, collecting \eqref{APRI_EQ_36}-\eqref{APRI_EQ_38} and using \eqref{APRI_EQ_28}, \eqref{APRI_EQ_35}, we see that
\begin{equation*}
\frac{\epsilon}{2}\norm{\grad\laplace\varphi}_{L^2(\L^2)}^2\leq C.
\end{equation*}
Together with \eqref{APRI_EQ_35} and using again elliptic regularity theory, this implies
\begin{equation}
\label{APRI_EQ_41}\norm{\varphi}_{L^2(H^3)}\leq C.
\end{equation}

\subsubsection{Regularity for the convection terms and the time derivatives}
By Hölder's and Young's inequalities and the continuous embedding $\H^1\hookrightarrow\L^6$, we observe that
\begin{equation*}
\norm{\grad\varphi\cdot\v}_{L^2(0,T;L^{\frac{3}{2}})}^2 \leq C\intT \norm{\v}_{\H^1}^2\norm{\grad\varphi}_{\L^2}^2\dt \leq C\norm{\grad\varphi}_{L^{\infty}(\L^2)}^2\norm{\v}_{L^2(0,T;\H^1)}^2.
\end{equation*}
Using \eqref{APRI_EQ_28}-\eqref{APRI_EQ_29} we see that

\begin{equation*}
\norm{\varphi\Gamma_{\v}}_{L^2(L^2)}^2 \leq C\norm{\varphi}_{L^{\infty}(H^1)}^2\norm{\Gamma_{\v}}_{L^2(L^3)}^2\leq C.
\end{equation*}
From the last two inequalities, we deduce that
\begin{equation}
\label{APRI_EQ_42}\norm{\divergence(\varphi\v)}_{L^2(0,T;L^{\frac{3}{2}})}\leq C.
\end{equation}
Finally, using \eqref{APRI_EQ_28}-\eqref{APRI_EQ_29} and \eqref{APRI_EQ_42}, a comparison argument in \eqref{WFORM_1b} shows that
\begin{equation}
\label{APRI_EQ_43}\norm{\delt\varphi}_{L^2((H^1)^*)}.
\end{equation}
Notice that we have lower time regularity for the time derivative of $\varphi$ compared to the convection term since the regularity of the time derivative depends on the term $\grad\mu$. Summarising all the estimates, we end up with
\begin{align}
\nonumber&\norm{\varphi}_{H^1((H^1)^*)\cap L^{\infty}(H^1)\cap L^4(H^2)\cap L^2(H^3) } + \norm{\sigma}_{L^2(H^1)} + \norm{\mu}_{L^2(H^1)}\\
\label{APRI_EQ_44}&\quad + \norm{\divergence(\varphi\v)}_{L^2(L^{\frac{3}{2}})}  + \norm{\v}_{L^2(\H^1)}+ \norm{p}_{L^2(L^2)}\leq C.
\end{align}

\subsection{Passing to the limit}
The a-priori-estimates \eqref{APRI_EQ_44} deduced within the Galerkin scheme are enough to deduce existence of solutions. We refer the reader to \cite{EbenbeckGarcke,GarckeLam1} for details in passing to the limit in the Galerkin scheme.

\subsection{Further results on regularity}
In the case when $\sigma_{\infty}\in L^4(L^2(\del\Omega))$, by \eqref{APRI_EQ_3} we obtain
\begin{equation}
\label{APRI_EQ_45}\norm{\sigma}_{L^4(H^1)}\leq  C.
\end{equation}
In particular, by \textnormal{(A4)} this gives
\begin{equation}
\label{APRI_EQ_46}\norm{\Gamma_{\v}}_{L^4(L^2)}+\norm{\Gamma_{\varphi}}_{L^4(L^2)}\leq C.
\end{equation}
Thanks to Lemma \ref{LEM_GAGNIR}, we have the continuous embedding
\begin{equation*}
L^{\infty}(H^1)\cap L^2(H^3)\hookrightarrow L^{20}(L^{10}).
\end{equation*}
Hence, the assumptions on $\psi(\cdot)$ and \eqref{APRI_EQ_44} imply 
\begin{equation}
\label{APRI_EQ_47}\norm{\psi'(\varphi)}_{L^4(L^2)}\leq C.
\end{equation}
Taking $\Phi = \mu+\chi\sigma$ in \eqref{WFORM_1c} and squaring the resulting identity, an application of Hölder's and Young's inequalities gives
\begin{equation*}
\norml{2}{\mu+\chi\sigma}^4 \leq C\left(\norml{2}{\psi'(\varphi)}^4 + \norm{\grad(\mu+\chi\sigma)}_{\L^2}^2\norm{\grad\varphi}_{\L^2}^2\right).
\end{equation*}
Integrating this inequality in time from $0$ to $T$ and using \eqref{APRI_EQ_44}, \eqref{APRI_EQ_47}, we conclude that
\begin{equation}
\label{APRI_EQ_48}\norm{\mu+\chi\sigma}_{L^4(L^2)}\leq C.
\end{equation}
We now choose $\boldsymbol{\Phi}=\v$ in \eqref{WFORM_1a}, use Young's, Hölder's and Korn's inequality (see \cite[Thm. 6.3-3]{Ciarlet}) together with the Sobolev embedding $\H^1\hookrightarrow \L^6$ to obtain
\begin{equation*}
\norm{\v}_{\H^1}^{\frac{8}{3}}\leq C\left(\norml{2}{p}^{\frac{4}{3}}\norml{2}{\Gamma_{\v}}^{\frac{4}{3}} + \norm{(\mu+\chi\sigma)\grad\varphi}_{\L^{\frac{6}{5}}}^{\frac{8}{3}}\right).
\end{equation*}
Integrating this inequality in time from $0$ to $T$, by Hölder's and Young's inequalities we get
\begin{equation}
\label{APRI_EQ_49}\norm{\v}_{L^{\frac{8}{3}}(\H^1)}^{\frac{8}{3}}\leq C\left(\norm{p}_{L^2(L^2)}^{\frac{4}{3}}\norm{\Gamma_{\v}}_{L^4(L^2)}^{\frac{4}{3}} + \norm{\mu+\chi\sigma}_{L^4(L^2)}^{\frac{8}{3}}\norm{\grad\varphi}_{L^8(\L^3)}^{\frac{8}{3}}\right).
\end{equation}
Applying Lemma \ref{LEM_GAGNIR}, we have the continuous embedding
\begin{equation}
\label{APRI_EQ_50}L^{\infty}(H^1)\cap L^4(H^2)\hookrightarrow L^8(W^{1,3}).
\end{equation}
Hence, using \eqref{APRI_EQ_44}, \eqref{APRI_EQ_47}  and  \eqref{APRI_EQ_48}, from \eqref{APRI_EQ_49} we conclude that
\begin{equation}
\label{APRI_EQ_51}\norm{\v}_{L^{\frac{8}{3}}(\H^1)}\leq C.
\end{equation}
Furthermore, using Hölder's and Young's inequalities, the continuous embedding $\H^1\hookrightarrow \L^6$ and \eqref{APRI_EQ_50}-\eqref{APRI_EQ_51} yields
\begin{equation*}
\norm{\grad\varphi\cdot\v}_{L^2(0,T;L^2)}^2 \leq C\intT \norm{\v}_{\H^1}^2\norm{\grad\varphi}_{\L^3}^2\dt \leq C\norm{\grad\varphi}_{L^8(\L^3)}^2\norm{\v}_{L^{\frac{8}{3}}(0,T;\H^1)}^2\leq C.
\end{equation*}
Together with the estimate (see \eqref{APRI_EQ_29}, \eqref{APRI_EQ_44})
\begin{equation*}
\norm{\varphi\Gamma_{\v}}_{L^2(L^2)}^2 \leq C\norm{\varphi}_{L^{\infty}(H^1)}^2\norm{\Gamma_{\v}}_{L^2(L^3)}^2\leq C,
\end{equation*}
this implies
\begin{equation}
\label{APRI_EQ_52}\norm{\divergence(\varphi\v)}_{L^2(L^2)}\leq C.
\end{equation}
Using \eqref{APRI_EQ_45}, \eqref{APRI_EQ_48}, \eqref{APRI_EQ_51}-\eqref{APRI_EQ_52} and recalling \eqref{APRI_EQ_44}, we obtain

\begin{align}
\nonumber&\norm{\varphi}_{H^1((H^1)^*)\cap L^{\infty}(H^1)\cap L^4(H^2)\cap L^2(H^3) } + \norm{\sigma}_{L^4(H^1)} + \norm{\mu}_{L^2(H^1)\cap L^4(L^2)}\\
\label{APRI_EQ_53}&\quad + \norm{\divergence(\varphi\v)}_{L^2(L^2)}  + \norm{\v}_{L^{\frac{8}{3}}(\H^1)}+ \norm{p}_{L^2(L^2)}\leq C.
\end{align}
\subsection{The singular limit of large boundary permeability} Due to Theorem \ref{THM_WSOL_1}, for every $K>0$ there exists a solution quintuple $(\varphi_K,\mu_K,\sigma_K,\v_K,p_K)$ solving \eqref{MEQ}, \eqref{BIC} in the sense of Definition \ref{DEF_WSOL_1} and enjoying the regularity properties stated in Theorem \ref{THM_WSOL_1}. In the following, we assume without loss of generality that $K>1$. Let $E\colon H^{\frac{1}{2}}(\del\Omega)\to H^1$ be a bounded, linear extension operator satisfying $(Ef)|_{\del\Omega} = f$ for all $f\in H^{\frac{1}{2}}(\del\Omega)$ (see \cite[Chap. 2, Thm. 5.7]{Necas} for details). Then, choosing $\Phi = \sigma_K-E\sigma_{\infty}$ in \eqref{WFORM_1d} (in the following, we will omit the operator $E$), we obtain
\begin{equation}
\label{limit_K_eq_1}\intO |\grad\sigma_K|^2 + h(\varphi)|\sigma_K|^2\dx  + K\int_{\del\Omega}|\sigma_K-\sigma_{\infty}|^2\d\mathcal{H}^{d-1} = \intO \grad\sigma_K\cdot\grad\sigma_{\infty}+h(\varphi)\sigma_K\sigma_{\infty} \dx .
\end{equation} 
For the first term on the r.h.s. of this equation, we use Hölder's and Young's inequalities and the boundedness of the extension operator to obtain
\begin{equation}
\label{limit_K_eq_2}\left|\intO \grad\sigma_K\cdot\grad\sigma_{\infty}\dx \right| \leq \frac{1}{4}\normL{2}{\grad\sigma_K}^2 + \normL{2}{\grad\sigma_{\infty}}^2\leq \frac{1}{4}\normL{2}{\grad\sigma_K}^2 + C\norm{\sigma_{\infty}}_{H^{\frac{1}{2}}(\del\Omega)}^2.
\end{equation}
With the same arguments and using the boundedness of $h(\cdot)$, we can estimate the second term on the r.h.s. of \eqref{limit_K_eq_1} by
\begin{equation}
\label{limit_K_eq_3}\left|\intO h(\varphi)\sigma_K\sigma_{\infty}\dx \right|\leq \delta\norml{2}{\sigma_K}^2 + C_{\delta}\norm{\sigma_{\infty}}_{H^{\frac{1}{2}}(\del\Omega)}^2.
\end{equation}
From Poincaré's inequality and the boundedness of the extension operator, we know that
\begin{equation}
\label{limit_K_eq_4}\norml{2}{\sigma_K}^2\leq \tilde{C}\left(\normL{2}{\grad\sigma_K}^2 + \norm{\sigma_K-\sigma_{\infty}}_{L^2(\del\Omega)}^2 + \norm{\sigma_{\infty}}_{H^{\frac{1}{2}}(\del\Omega)}^2\right),
\end{equation}
for a positive constant $\tilde{C}$ independent of $K$.
Choosing $\delta$ small enough and using this inequality in \eqref{limit_K_eq_3}, we obtain 
\begin{equation}
\label{limit_K_eq_5}\left|\intO h(\varphi)\sigma_K\sigma_{\infty}\dx \right|\leq  \frac{1}{4}\left(\normL{2}{\grad\sigma_K}^2 + \norm{\sigma_K-\sigma_{\infty}}_{L^2(\del\Omega)}^2\right) + C\norm{\sigma_{\infty}}_{H^{\frac{1}{2}}(\del\Omega)}^2.
\end{equation}
Plugging in \eqref{limit_K_eq_2}, \eqref{limit_K_eq_5} into \eqref{limit_K_eq_1} and neglecting the non-negative term $\intO h(\varphi)|\sigma_K|^2\dx $ on the l.h.s. of \eqref{limit_K_eq_1}, we arrive at
\begin{equation}
\label{limit_K_eq_6}\intO |\grad\sigma_K|^2 \dx  + K\int_{\del\Omega}|\sigma_K-\sigma_{\infty}|^2\d\mathcal{H}^{d-1}\leq C\norm{\sigma_{\infty}}_{H^{\frac{1}{2}}(\del\Omega)}^2.
\end{equation}
Multiplying \eqref{limit_K_eq_4} by $\frac{1}{2\tilde{C}}$ and adding the resulting equation to \eqref{limit_K_eq_6} yields
\begin{equation*}
\normh{1}{\sigma_K}^2 + K\int_{\del\Omega}|\sigma_K-\sigma_{\infty}|^2\d\mathcal{H}^{d-1}\leq C\norm{\sigma_{\infty}}_{H^{\frac{1}{2}}(\del\Omega)}^2.
\end{equation*}
Integrating this inequality in time from $0$ to $T$ and using $\sigma_{\infty}\in L^2(H^{\frac{1}{2}}(\del\Omega))$, we conclude that
\begin{equation*}
\norm{\sigma_K}_{L^2(H^1)} + \sqrt{K}\norm{\sigma_K-\sigma_{\infty}}_{L^2(L^2(\del\Omega))}\leq C,
\end{equation*}
where $C$ is independent of $K$. Then, with exactly the same arguments as above, it follows that
\begin{align}
\nonumber&\norm{\varphi_K}_{H^1((H^1)^*)\cap L^{\infty}(H^1)\cap L^4(H^2)\cap L^2(H^3) } + \norm{\sigma_K}_{L^2(H^1)}+ \sqrt{K}\norm{\sigma_K-\sigma_{\infty}}_{L^2(L^2(\del\Omega))}\\
\label{limit_K_eq_7}&\quad   + \norm{\mu_K}_{L^2(H^1)}+\norm{\divergence(\varphi_K\v_K)}_{L^2(L^{\frac{3}{2}})}+\norm{\v_K}_{L^{2}(\H^1)}+ \norm{p_K}_{L^2(L^2)}\leq C.
\end{align}
Using  standard compactness arguments (Aubin-Lions theorem (see \cite[Sec. 8, Cor. 4]{Simon}) and reflexive weak compactness), we obtain exactly the convergence results as stated in Theorem \ref{THM_KLIM}. Passing to the limit can be carried out with exactly the same arguments as stated in Section 3.2. We will only present the arguments needed for \eqref{WFORM_1d}. In the following, let $\xi\in H_0^1$ be arbitrary. Multiplying \eqref{WFORM_1d} with $\delta\in C_0^{\infty}(0,T)$, integrating in time from $0$ to $T$ and noting that $H_0^1\subset H^1$, we observe that
\begin{equation}
\label{limit_K_eq_8}0 = \intT\intO \delta (\grad\sigma_K\cdot\grad\xi + h(\varphi_K)\sigma_K\xi)\dx \dt \quad\forall \xi\in H_0^1.
\end{equation}
Since $h(\cdot)$ is a bounded, continuous function, $\delta\xi\in C^{\infty}(H_0^1)$ and $\varphi_K\to\varphi$ a.\,e. in $Q$, the Lebesgue theorem gives that
\begin{equation*}
\norm{h(\varphi_K)\delta\xi-h(\varphi)\delta\xi}_{L^2(Q)}\to 0\quad\text{as }K\to\infty.
\end{equation*}
Since $\sigma_K\to\sigma$ weakly in $L^2(Q)$ as $K\to\infty$, by the product of weak-strong convergence we obtain
\begin{equation}
\label{limit_K_eq_9}\intT\intO h(\varphi_K)\sigma_K\xi\dx \dt \to \intT\intO h(\varphi)\sigma\xi\dx \dt \quad\text{as }K\to\infty.
\end{equation}
Furthermore, since $\sigma_K\to\sigma$ weakly in $L^2(H^1)$ and as $\delta\xi\in L^2(H^1)$, it follows that
\begin{equation}
\label{limit_K_eq_10}\intT\intO \delta\grad\sigma_K\cdot\grad\xi\dx \dt \to \intT\intO \delta\grad\sigma\cdot\grad\xi\dx \dt \quad\text{as }K\to\infty.
\end{equation}
Due to \eqref{limit_K_eq_9}-\eqref{limit_K_eq_10}, we can pass to the limit in \eqref{limit_K_eq_8} to deduce that
\begin{equation*}
0 = \intT\intO \delta (\grad\sigma\cdot\grad\xi + h(\varphi)\sigma\xi)\dx \dt \quad\forall \xi\in H_0^1.
\end{equation*}
Since this holds for all $\delta\in C_0^{\infty}(0,T)$, we can recover \eqref{WF_KLIM_4}. Finally, from \eqref{limit_K_eq_7}, we know that
\begin{equation*}
\norm{\sigma_K-\sigma_{\infty}}_{L^2(L^2(\del\Omega))}\leq \frac{C}{\sqrt{K}},
\end{equation*}
where $C$ is independent of $K$. Letting $K\to\infty$ and recalling that $\sigma_K\to\sigma$ weakly in $L^2(L^2(\del\Omega))$ as $K\to\infty$, it follows that
\begin{equation*}
\sigma = \sigma_{\infty}\quad\text{a.\,e. on }\Sigma,
\end{equation*}
which completes the proof.
\section{The singular limit of vanishing viscosities}
Let $\{\eta_n,\lambda_n\}_{n\in\N}$ be a sequence of function pairs fulfilling \textnormal{(A3)} such that
\begin{equation}
\label{VISC_LIM_1}\norm{\eta_n(\cdot)}_{C^0(\R)}\to 0,\quad \norm{\lambda_n(\cdot)}_{C^0(\R)}\to 0\quad\text{as }n\to\infty.
\end{equation}
Without loss of generality, we can assume that
\begin{equation}
\label{VISC_LIM_2}\norm{\eta_{n}(\cdot)}_{L^{\infty}(\R)}\leq 1,\quad \norm{\lambda_{n}(\cdot)}_{L^{\infty}(\R)}\leq 1.
\end{equation}
Then, by Theorem \ref{THM_WSOL_1}, for every $n\in\N$ there exists a solution quintuple $(\varphi_{n},\mu_{n},\sigma_{n},\v_{n},p_{n})$ of \eqref{MEQ}, \eqref{BIC} in the sense of Definition \ref{DEF_WSOL_1} fulfilling
\begin{align*}
&\varphi_{n} \in  H^1(0,T;(H^1)^*)\cap L^{\infty}(0,T;H^1)\cap L^4(0,T;H^2)\cap L^2(0,T;H^3),\quad \mu_{n}\in L^2(0,T;H^1),\\
&\sigma_{n}\in L^2(0,T;H^1),\quad \v_{n}\in L^2(0,T;\H^1),\quad p_{n}\in L^{2}(0,T;L^2),
\end{align*}
such that 
\begin{equation}
\label{VISC_LIM_3}\divergence(\v_{n}) = \Gamma_{\v}(\varphi_n,\sigma_n)\text{ a.\,e. in }Q,\quad \deln\varphi_{n} = 0\text{ a.\,e. in }\Sigma,\quad \varphi_{n}(0)=\varphi_0\text{ a.\,e. in }\Omega,
\end{equation}
and 
\begin{subequations}
	\label{VISC_LIM_4}
	\begin{align}
	\label{VISC_LIM_4a}0&= \intO \T_n(\v_{n},p_{n})\colon \grad\boldsymbol{\Phi}+\nu\v_{n}\cdot\boldsymbol{\Phi} -  (\mu_{n}+\chi\sigma_{n})\grad\varphi_{n}\cdot\boldsymbol{\Phi}\d x,\\
	\label{VISC_LIM_4b} 0&= \langle\delt\varphi_{n}{,}\Phi\rangle_{H^1,(H^1)^*} + \intO m(\varphi_{n})\grad\mu_{n}\cdot\grad\Phi + (\grad\varphi_{n}\cdot\v_n  +\varphi_n\Gamma_{\v}(\varphi_n,\sigma_n)-\Gamma_{\varphi}(\varphi_n,\sigma_n))\Phi\d x,\\
	\label{VISC_LIM_4c} 0  &= \intO \grad\sigma_{n}\cdot\grad\Phi + h(\varphi_{n})\sigma_{n}\Phi\d x + \int_{\del\Omega}K(\sigma_{n}-\sigma_{\infty})\Phi\d\mathcal{H}^{d-1},
	\end{align}
for a.\,e. $t\in(0,T)$ and for all $\boldsymbol{\Phi}\in \H^1,\,\Phi\in H^1$, where $\mu_n$ is given by
\begin{equation}
\label{VISC_LIM_4d}\mu_{n} =  \epsilon^{-1}\Psi'(\varphi_{n})- \epsilon\laplace\varphi_{n} - \chi\sigma_{n}\quad \text{a.\,e in }Q
\end{equation}
\end{subequations}
and 
\begin{equation*}
\T_n(\v_{n},p_{n})\coloneqq 2\eta_{n}(\varphi_{n})  \D\v_{n} + \lambda_{n}(\varphi_{n})\divergence(\v_{n})\mathbf{I} - p_{n}\I.
\end{equation*}
We will denote $\Gamma_{\v,n} =\Gamma_{\v}(\varphi_n,\sigma_n),\,\Gamma_{\varphi,n} =\Gamma_{\varphi}(\varphi_n,\sigma_n)$. 
\subsection{A-priori-estimates}
In the following, we will derive estimates which are independent of $n\in\N$. By $C$, we will denote a generic constant depending on the system parameters and on $\Omega$, $T$, but not on $n\in\N$. Furthermore, we will frequently use Hölder's and Young's inequalities.\newline 
First, we recall that \textnormal{(A1)}, \textnormal{(A5)} and the continuous embedding $H^1\hookrightarrow L^6$ imply that $\psi(\varphi_0)\in L^1,\, \nabla\varphi_0\in \L^2$. 
Then, using \textnormal{(A1)}, \eqref{APRI_EQ_25a}-\eqref{APRI_EQ_26} and \eqref{VISC_LIM_2}, taking the supremum over all $s\in (0,T]$ in \eqref{APRI_EQ_26} yields
\begin{align}
\nonumber &  \esssup_{s\in (0,T]}\left(\norml{1}{\psi(\varphi_n(s))} + \normL{2}{\grad\varphi_n(s)}^2\right) \\
\label{VISC_LIM_5}&\quad + \int_{0}^{T}\frac{m_0}{4}\normL{2}{\grad\mu_n}^2+  \frac{\nu}{2}\normL{2}{\v_n}^2 + \normL{2}{\sqrt{\eta_n(\varphi_n)}\D\v_n}^2 + \normh{1}{\sigma_n}^2\d t\leq C.
\end{align}
Recalling \eqref{ASS_PSI_4}, using Poincaré's inequality, \eqref{APRI_EQ_12}, \eqref{VISC_LIM_5} and the fact that $\rho\geq 2$, this in particular gives
\begin{equation}
\label{VISC_LIM_6}\esssup_{s\in (0,T]}\normh{1}{\varphi_n(s)}^2 + \intT \normh{1}{\mu_n}^2\dt\leq C.
\end{equation}
Now, using exactly the same arguments as in Subsection 3.1.6, we obtain
\begin{equation}
\label{VISC_LIM_7}\norm{\varphi_n}_{L^4(H^2)}\leq C.
\end{equation}
Then, using \textnormal{(A5)}, it is straightforward to check that $\psi'(\varphi)\in L^2(H^1)$ with bounded norm. Together with \eqref{VISC_LIM_5}-\eqref{VISC_LIM_7} and using elliptic regularity theory, this implies
\begin{equation}
\label{VISC_LIM_7a}\norm{\varphi_n}_{L^4(H^2)\cap L^2(H^3)}\leq C.
\end{equation}
By the specific form of $\Gamma_{\v}(\cdot{,}\cdot),\,\Gamma_{\varphi}(\cdot{,}\cdot)$ and the continuous embedding $H^1\hookrightarrow L^6$, an application of \eqref{VISC_LIM_3} and \eqref{VISC_LIM_5} yields
\begin{equation}
\label{VISC_LIM_8}\norm{\Gamma_{\v,n}}_{L^2(L^6)} + \norm{\Gamma_{\varphi,n}}_{L^2(L^6)}\leq C.
\end{equation}
Using \eqref{APRI_EQ_32} and \eqref{VISC_LIM_2}, for every $n\in\N$ we obtain
\begin{equation*}
\norml{2}{p_n}^2 \leq C\left(\normL{2}{\sqrt{\eta_n(\varphi_n)}\D\v_n}^2 + \norml{2}{\Gamma_{\v,n}}^2 + \normL{2}{\v_n}^2 + \norml{3}{\mu_n + \chi\sigma_n}^2\normL{2}{\grad\varphi_n}^2\right).
\end{equation*}
Integrating this inequality in time from $0$ to $T$ and using \eqref{VISC_LIM_5}-\eqref{VISC_LIM_6}, \eqref{VISC_LIM_8}, we have
\begin{equation}
\label{VISC_LIM_9}\norm{p_n}_{L^2(L^2)}\leq C.
\end{equation}
Using \eqref{VISC_LIM_6}-\eqref{VISC_LIM_8} and the continuous embedding $H^1\hookrightarrow L^3$, we calculate
\begin{equation}
\label{VISC_LIM_10}\norm{\varphi_n\Gamma_{\v,n}}_{L^2(L^2)}\leq \norm{\varphi_n}_{L^{\infty}(L^3)}\norm{\Gamma_{\v,n}}_{L^2(L^6)}\leq C\norm{\varphi_n}_{L^{\infty}(H^1)}\norm{\Gamma_{\v,n}}_{L^2(L^6)}\leq C.
\end{equation} 
Now, let $\rho\in L^{\frac{8}{3}}(H^1)$. Then, using \eqref{LEM_GAGNIR_EST} and \eqref{VISC_LIM_5}-\eqref{VISC_LIM_7a}, we obtain
\begin{align*}
\left|\intT\intO \nabla\varphi_n\cdot\v_n\rho\dx \dt\right| &\leq \intT \normL{3}{\nabla\varphi_n}\normL{2}{\v_n}\normh{1}{\rho}\dt \\
&\leq \intT \normh{1}{\varphi_n}^{\frac{3}{4}}\normh{3}{\varphi_n}^{\frac{1}{4}}\normL{2}{\v_n}\normh{1}{\rho}\dt \\
&\leq \norm{\varphi_n}_{L^{\infty}(H^1)}^{\frac{3}{4}}\norm{\varphi_n}_{L^2(H^3)}^{\frac{1}{4}}\norm{\v_n}_{L^2(\L^2)}\norm{\rho}_{L^{\frac{8}{3}}(H^1)}\\
&\leq C\norm{\rho}_{L^{\frac{8}{3}}(H^1)},
\end{align*}
which implies
\begin{equation*}
\norm{\nabla\varphi_n\cdot\v_n}_{L^{\frac{8}{5}}((H^1)^*)}\leq C.
\end{equation*}
Using the continuous embedding $L^2\hookrightarrow (H^1)^*$ and \eqref{VISC_LIM_3}, \eqref{VISC_LIM_10}, we deduce that
\begin{equation}
\label{VISC_LIM_11}\norm{\divergence(\varphi_n\v_n)}_{L^{\frac{8}{5}}((H^1)^*)}\leq C.
\end{equation}
Finally, a comparison argument in \eqref{VISC_LIM_4b} yields
\begin{equation}
\label{VISC_LIM_12}\norm{\delt\varphi_n}_{L^{\frac{8}{5}}((H^1)^*)}\leq C,
\end{equation}
where we used \eqref{VISC_LIM_5}-\eqref{VISC_LIM_8} and \eqref{VISC_LIM_11}. Summarising \eqref{VISC_LIM_5}-\eqref{VISC_LIM_12}, we end up with
\begin{align}
\nonumber &\norm{\varphi_n}_{W^{1,\frac{8}{5}}((H^1)^*)\cap L^{\infty}(H^1)\cap L^4(H^2)\cap L^2(H^3)} + \norm{\mu_n}_{L^2(H^1)} + \norm{\sigma_n}_{L^2(H^1)} + \norm{\divergence(\varphi_n\v_n)}_{L^{\frac{8}{5}}((H^1)^*)}\\
\label{VISC_LIM_13}&\quad + \norm{\divergence(\v_n)}_{L^2(L^6)}+ \norm{\sqrt{\eta_n(\varphi_n)}\D\v_n}_{L^2(\L^2)} + \norm{\v_n}_{L^2(\L_{\text{\divergence}}^2(\Omega))} + \norm{p_n}_{L^2(L^2)}\leq C.
\end{align}

\subsection{Passing to the limit}
Recalling \eqref{TRACE_ESTIM} and \eqref{VISC_LIM_13}, using  standard compactness arguments (Aubin-Lions theorem (see \cite[Sec. 8, Cor. 4]{Simon}) and reflexive weak compactness), the compact embeddings
\begin{equation*}
H^{j+1}(\Omega) = W^{j+1,2}(\Omega)\hookrightarrow\hookrightarrow W^{j,r}\quad \forall j\in \mathbb{Z},\,j\geq 0,\,1\leq r<6,
\end{equation*}
and $L^2\hookrightarrow\hookrightarrow (H^1)^*$, we obtain, at least for a subsequence which will again be labelled by $n$, the following convergence results:
\begin{alignat*}{3}
\varphi_{n}&\to\varphi&&\quad\text{weakly-}*&&\quad\text{ in } W^{1,\frac{8}{5}}((H^1)^*)\cap L^{\infty}(H^1)\cap L^4(H^2)\cap  L^2(H^3),\\
\sigma_{n}&\to \sigma &&\quad\text{weakly}&&\quad\text{ in }L^2(H^1),\\
\mu_{n}&\to \mu&&\quad\text{weakly}&&\quad\text{ in }L^2(H^1),\\
p_{n}&\to p&&\quad\text{weakly}&&\quad\text{ in }L^{2}(L^2),\\
\v_{n}&\to\v&&\quad\text{weakly}&&\quad\text{ in }L^2(\L^2)\cap L^2(L_{\divergence}^2(\Omega)),\\
\divergence(\v_{n})&\to\xi &&\quad\text{weakly}&&\quad\text{ in }L^2(L^2),\\
\v_{n}\cdot\n&\to\v\cdot\n&&\quad\text{weakly}&&\quad\text{ in }L^2((H^{\frac{1}{2}}(\del\Omega))^*),\\
\divergence(\varphi_{n}\v_{n})&\to \tau&&\quad\text{weakly}&&\quad\text{ in }L^{\frac{8}{5}}((H^1)^*),
\end{alignat*}
for some limit functions $\xi\in L^2(L^2),\,\tau\in L^{\frac{8}{5}}((H^1)^*)$. Furthermore,  by (\ref{VISC_LIM_1}), we have the strong convergences
\begin{alignat*}{3}
&\varphi_{n}\to\varphi\quad\text{strongly}&&\quad\text{ in }C^0(L^r)\cap L^2(W^{2,r})\text{ and a.e. in }Q,\\
&\eta_n(\varphi_n)\to 0&&\quad\text{ in }Q,\\
&\lambda_n(\varphi_n)\to 0&&\quad\text{ in }Q,
\end{alignat*}
for $r\in [1,6)$. In the following, we fix $\delta\in C_0^{\infty}(0,T),\,\Phi\in H^1,\,\boldsymbol{\Phi}\in \H^1$ and we note that
$\delta\Phi\in C^{\infty}(H^1),\,\delta\boldsymbol{\Phi}\in C^{\infty}(\H^1)$. Multiplying \eqref{VISC_LIM_4a}-\eqref{VISC_LIM_4d} with $\delta$  and integrating in time from $0$ to $T$ and multiplying \eqref{VISC_LIM_4a} with $\delta\Phi$ and integrating over $Q$, we obtain
\begin{subequations} \label{VISC_LIM_14}
	\begin{align}
	\nonumber 0 &= \intT\intO \delta(t)(2\eta_n(\varphi_n) \D\v_{n} + \lambda_n(\varphi_n)\divergence(\v_n)\I - p_{n}\I)\colon \grad\boldsymbol{\Phi}+\nu\v_{n}\cdot\boldsymbol{\Phi})\d x\d t \\
	\label{VISC_LIM_14a}&- \intT\intO \delta(t)(\mu_{n}+\chi\sigma_{n})\grad\varphi_{n}\cdot\boldsymbol{\Phi}\d x\d t,\\
	\nonumber 0& =\intT \delta(t)\langle\delt\varphi_{n}{,}\Phi\rangle_{H^1,(H^1)^*}\d t  + \intT \intO \delta(t)(m(\varphi_{n})\grad\mu_{n}\cdot\grad\Phi - \Gamma_{\varphi,n}\Phi)\d x\d t\\
	\label{VISC_LIM_14b} &\quad +\intT \intO \delta(t)(\grad\varphi_n\cdot\v_n  +\varphi_n\Gamma_{\v,n})\Phi\d x\d t,\\
	\label{VISC_LIM_14c}0&=\intT\intO \delta(t)\mu_{n}\Phi\d x\d t -  \intT\intO \delta(t)(\epsilon^{-1}\Psi'(\varphi_{n})- \epsilon\laplace\varphi_{n} - \chi\sigma_{n})\Phi\d x\d t,\\
	\label{VISC_LIM_14d} 0&   = \intT\intO \delta(t)(\grad\sigma_{n}\cdot\grad\Phi + h(\varphi_{n})\sigma_{n}\Phi)\d x\d t + \intT \int_{\del\Omega}\delta(t)(K(\sigma_{n}-\sigma_{\infty})\Phi)\d\mathcal{H}^{d-1}\d t.
	\end{align}
	Furthermore, we multiply $\eqref{VISC_LIM_3}_1$ with $\delta\Phi$ and integrate over $Q$ to obtain
	\begin{equation}
	\label{VISC_LIM_14e}\intT\intO \delta \divergence(\v_{n})\Phi\d x\d t = \intT\intO \delta(t)\Gamma_{\v,n}\Phi\d x\d t.
	\end{equation}
\end{subequations}
We now want to analyse each term individually. For \eqref{VISC_LIM_14c}-\eqref{VISC_LIM_14d}, we omit the details and refer to the arguments used in \cite[Sec. 5]{EbenbeckGarcke}, \cite[Sec. 5]{GarckeLam1}. 
\paragraph{Step 1 (\eqref{VISC_LIM_14e}):}
 Since $\varphi_{n}\to \varphi$ a.\,e. in $Q$ and due to the boundedness of $b_{\v}(\cdot),\, f_{\v}(\cdot)$, Lebesgue dominated convergence theorem implies
\begin{equation*}
\norm{\delta\Phi(b_{\v}(\varphi_{n})-b_{\v}(\varphi))}_{L^2(Q)}\to 0,\quad \norm{\delta\Phi(f_{\v}(\varphi_{n})-f_{\v}(\varphi))}_{L^2(Q)}\to 0
\end{equation*}
as $n\to\infty$.
Together with the weak convergence $\sigma_{n}\to \sigma$ in $L^2(Q)$ as $n\to\infty$, by the product of weak-strong convergence we obtain
\begin{equation}
\label{VISC_LIM_15} \intT\intO \delta(t)\Gamma_{\v,n}\Phi\d x\d t\to \intT\intO \delta(t)\Gamma_{\v}(\varphi,\sigma)\Phi\d x\d t\quad\text{as }n\to\infty.
\end{equation}
Using \eqref{TRACE_GENERAL}, we see that
\begin{equation*}
\intT\intO \delta\divergence(\v_{n})\Phi\d x\d t = \intT\delta(t) \langle \v_{n}\cdot\n{,}\Phi\rangle_{H^{\frac{1}{2}}(\del\Omega)}\d t - \intT\intO \delta(t)\v_{n}\cdot\grad\Phi\d x\d t.
\end{equation*}
Since $\v_{n}\cdot\n\to\v\cdot\n$ weakly in $L^2((H^{\frac{1}{2}}(\del\Omega))^*)$, $\v_{n}\to\v$ weakly in $L^2(\L^2)$ and $\divergence(\v_{n})\to\xi$ weakly in $L^2(L^2)$ as $n\to\infty$, we can pass to the limit on both sides of this equation to obtain
\begin{equation*}
\intT\intO \delta\xi \Phi\d x\d t = \intT\delta(t) \langle \v\cdot\n{,}\Phi\rangle_{H^{\frac{1}{2}}(\del\Omega)}\d t - \intT\intO \delta(t)\v\cdot\grad\Phi\d x\d t.
\end{equation*}
Since $\v\in L_{\divergence}^2(\Omega)$, we can again use \eqref{TRACE_GENERAL} to obtain
\begin{equation*}
\intT\intO \delta\xi \Phi\d x\d t = \intT \intO \delta\divergence(\v)\Phi\d x\d t.
\end{equation*}
This in particular gives
\begin{equation}
\label{VISC_LIM_16}\divergence(\v) = \xi \quad \text{a.e. in }Q.
\end{equation}
From this considerations and recalling \eqref{VISC_LIM_15}, we can pass to the limit $n\to\infty$ in \eqref{VISC_LIM_14e} to obtain
\begin{equation*}
\intT\intO \delta\divergence(\v)\Phi\d x\d t = \intT\intO \delta \Gamma_{\v}(\varphi,\sigma)\Phi\d x\d t,
\end{equation*}
which in particular implies 
\begin{equation}
\label{VISC_LIM_17}\divergence(\v) = \Gamma_{\v}(\varphi,\sigma)\quad\text{a.\,e. in }Q.
\end{equation}
\paragraph{Step 2 (\eqref{VISC_LIM_14b}):} Since $\delta\Phi\in C^{\infty}(H^1)$ and $\divergence(\varphi_{n}\v_{n})\to \tau$ weakly in $L^{\frac{8}{5}}((H^1)^*)$, we have
\begin{equation}
\label{VISC_LIM_18}\intT \intO\delta\divergence(\varphi_{n}\v_{n})\Phi\d x\d t \to \intT \delta(t)\langle\tau{,}\Phi\rangle_{H^1}\d t\quad \text{as }n\to\infty.
\end{equation}
Moreover, since $\grad\varphi_{n}\to\grad\varphi$ strongly in $L^2(\L^3)$ and due to the continuous embedding $H^1\hookrightarrow L^6$, we have
\begin{align*}
\intT\intO |\delta|^2|\Phi|^2|\grad\varphi_{n}-\grad\varphi|^2\d x\d t&\leq \intT |\delta|^2 \norml{6}{\Phi}^2\normL{3}{\grad\varphi_{n}-\grad\varphi}^2\d t\\
&\leq C\norm{\delta}_{L^{\infty}(0,T)}^2\normh{1}{\Phi}^2\norm{\grad\varphi_{n}-\grad\varphi}_{L^2(\L^3)}^2\to 0
\end{align*}
as $n\to\infty$. This implies $\delta\Phi\grad\varphi_{n}\to\delta\phi\grad\varphi$ strongly in $L^2(\L^2)$. Together with the weak convergence $\v_{n}\to\v$ in $L^2(\L^2)$ as $n\to\infty$, by the product of weak-strong convergence we get
\begin{equation}
\label{VISC_LIM_19}\intT\intO \delta\Phi\grad\varphi_{n}\cdot\v_{n}\d x\d t\to \intT\intO \delta\Phi\grad\varphi\cdot\v\d x\d t\quad\text{as }n\to\infty.
\end{equation}
Since $\varphi_{n}\to\varphi$ strongly in $L^2(L^3)$ and a.\,e. in $Q$ as $n\to\infty$, the boundedness of $b_{\v}(\cdot),\,f_{\v}(\cdot)$ and Lebesgue dominated convergence theorem imply
\begin{equation*}
\norm{(b_{\v}(\varphi_{n})\varphi_{n}-b_{\v}(\varphi)\varphi)\delta\Phi}_{L^2(Q)}\to 0,\qquad \norm{(f_{\v}(\varphi_{n})\varphi_{n}-f_{\v}(\varphi)\varphi)\delta\Phi}_{L^2(Q)}\to 0
\end{equation*}
as $n\to \infty$, where we used that $\Phi\in H^1\hookrightarrow L^6$. Together with the weak convergence $\sigma_{n}\to\sigma$ in $L^2(Q)$ as $n\to\infty$, this implies
\begin{equation}
\label{VISC_LIM_20}\intT\intO \delta\Phi \Gamma_{\v,n}\varphi_{n}\d x\d t\to \intT\intO \delta\Phi \Gamma_{\v}(\varphi,\sigma)\varphi\d x\d t\quad  \text{as }n\to\infty.
\end{equation}
Using $\eqref{VISC_LIM_3}_1$, we see that
\begin{equation*}
\intT \intO \delta\divergence(\varphi_{n}\v_{n})\Phi\d x\d t = \intT\intO \delta\Phi\grad\varphi_{n}\cdot\v_{n}\d x\d t + \intT\intO \delta\Phi \Gamma_{\v,n}\varphi_{n}\d x\d t.
\end{equation*}
Passing to the limit $n\to\infty$ on both sides of this equation and using \eqref{VISC_LIM_19}-\eqref{VISC_LIM_20}, we obtain
\begin{equation}
\label{VISC_LIM_21}\intT \delta(t)\langle\tau{,}\Phi\rangle_{H^1}\d t = \intT\intO \delta\Phi\grad\varphi\cdot\v\d x\d t + \intT\intO \delta\Phi \Gamma_{\v}(\varphi,\sigma)\varphi\d x\d t.
\end{equation}
Together with \eqref{VISC_LIM_17}, this gives
\begin{equation}
\label{VISC_LIM_22}\intT \delta(t)\langle\tau{,}\Phi\rangle_{H^1,(H^1)^*}\d t = \intT\intO \delta(t)\divergence(\varphi\v)\Phi\d x\d t,
\end{equation}
hence $\divergence(\varphi\v) = \tau$ in the sense of distributions.
For the remaining terms in \eqref{VISC_LIM_14b}, we again refer to \cite[Sec. 5]{EbenbeckGarcke}, \cite[Sec. 5]{GarckeLam1}, where they used similar arguments.
\paragraph{Step 3 (\eqref{VISC_LIM_14a}):} With exactly the same arguments as used for \eqref{VISC_LIM_19}, it follows that $\delta\boldsymbol{\Phi}\cdot\grad\varphi_{n}\to \delta\boldsymbol{\Phi}\cdot\grad\varphi$ strongly in $L^2(L^2)$ as $n\to\infty$. Then, recalling that $\mu_{n} +\chi\sigma_{n}\to \mu + \chi\sigma$ weakly in $L^2(L^2)$ as $n\to\infty$, by the product of weak-strong convergence we obtain
\begin{equation}
\label{VISC_LIM_23}\intT\intO \delta(t)(\mu_{n}+\chi\sigma_{n})\grad\varphi_{n}\cdot\boldsymbol{\Phi}\d x\d t\to \intT\intO \delta(t)(\mu+\chi\sigma)\grad\varphi\cdot\boldsymbol{\Phi}\d x\d t\quad\text{as }n\to\infty.
\end{equation}
Now, since $p_{n}\to p,\,\v_{n}\to\v$ weakly in $L^2(L^2)$ and $L^2(\L^2)$ respectively and observing
\begin{equation*}
\intT\intO \delta p_{n}\I\colon\grad\boldsymbol{\Phi}\d x\d t = \intT\intO \delta p_{n}\divergence(\boldsymbol{\Phi})\d x\d t,
\end{equation*}
we obtain that
\begin{equation}
\label{VISC_LIM_24}\intT\intO \delta (-p_{n}\divergence(\boldsymbol{\Phi}) +\nu\v_{n}\cdot\boldsymbol{\Phi})\d x\d t\to \intT\intO \delta (-p\divergence(\boldsymbol{\Phi}) +\nu\v\cdot\boldsymbol{\Phi})\d x\d t\quad\text{as }n\to \infty.
\end{equation}
Finally, we recall that $\eta_n(\varphi_n) \to 0$ a.\,e. in $Q$ as $n\to\infty$. Consequently, applying \eqref{VISC_LIM_13} yields
\begin{align}
\nonumber \left|\intT\intO \delta(t)2\eta_n(\varphi_n) \D\v_{n}\colon\grad\boldsymbol{\Phi}\d x\d t\right| &\leq C\norm{\sqrt{\eta_n(\varphi_n)}\D\v_n}_{L^2(\L^2)}\norm{\sqrt{\eta_n(\varphi_n)}}_{L^{\infty}(Q)}\norm{\delta}_{L^{\infty}(0,T)}\normH{1}{\boldsymbol{\Phi}}\\
\nonumber &\leq C \norm{\sqrt{\eta_n(\varphi_n)}}_{L^{\infty}(Q)}\norm{\delta}_{L^{\infty}(0,T)}\normH{1}{\boldsymbol{\Phi}}\\
\label{VISC_LIM_25} &\to 0\quad\text{as }n \to \infty.
\end{align}
Using that $\lambda_n(\varphi_n) \to 0$ a.\,e. in $Q$ as $n\to\infty$ and applying \eqref{VISC_LIM_13}, it follows that
\begin{align}
\nonumber \left|\intT\intO \delta(t)\lambda_n(\varphi_n) \divergence(\v_n)\I\colon\grad\boldsymbol{\Phi}\d x\d t\right| &\leq C\norm{\divergence(\v_n)}_{L^2(L^2)}\norm{\lambda_n(\varphi_n)}_{L^{\infty}(Q)}\norm{\delta}_{L^{\infty}(0,T)}\normH{1}{\boldsymbol{\Phi}}\\
\nonumber &\leq C \norm{\lambda_n(\varphi_n)}_{L^{\infty}(Q)}\norm{\delta}_{L^{\infty}(0,T)}\normH{1}{\boldsymbol{\Phi}}\\
\label{VISC_LIM_26} &\to 0\quad\text{as }n \to \infty.
\end{align}
\paragraph{Step 4:} Due to \eqref{VISC_LIM_15}-\eqref{VISC_LIM_26}, we have enough to pass to the limit $n\to\infty$ in \eqref{VISC_LIM_14} to obtain that
\begin{subequations}\label{VISC_LIM_27}
\begin{align}
\label{VISC_LIM_27a}0&=\intT\delta(t)\left(\intO -p\divergence(\boldsymbol{\Phi}) +(\nu\v-(\mu+\chi\sigma)\grad\varphi)\cdot\boldsymbol{\Phi}\d x\right)\d t,\\
\label{VISC_LIM_27b} 0&=\intT \delta(t)\left(\langle\delt\varphi{,}\Phi\rangle_{H^1} + \intO m(\varphi)\grad\mu\cdot\grad\Phi + (\grad\varphi\cdot\v  +\varphi\Gamma_{\v}(\varphi,\sigma)-\Gamma_{\varphi}(\varphi,\sigma))\Phi\d x\right)\d t,\\
\label{VISC_LIM_27c}0&=\intT\delta(t)\left(\intO (\mu-\epsilon^{-1}\Psi'(\varphi)+ \epsilon\laplace\varphi + \chi\sigma)\Phi\d x\right)\d t,\\
\label{VISC_LIM_27d} 0&   = \intT \delta(t)\left(\intO \grad\sigma\cdot\grad\Phi + h(\varphi)\sigma\Phi\d x + \int_{\del\Omega}K(\sigma-\sigma_{\infty})\Phi\d\mathcal{H}^{d-1}\right)\d t,
\end{align}
\end{subequations}
for all $\delta\in C_0^{\infty}(0,T)$ and
\begin{subequations} \label{VISC_LIM_28}
\begin{equation}
\label{VISC_LIM_28a}\divergence(\v)=\Gamma_{\v}(\varphi,\sigma)\quad \text{a.\,e. in }Q.
\end{equation}
Since (\ref{VISC_LIM_27}) holds for all $\delta\in C_0^{\infty}(0,T)$, we deduce that 
\begin{align}
\label{VISC_LIM_28b}0&=\intO -p\divergence(\boldsymbol{\Phi}) +(\nu\v-(\mu+\chi\sigma)\grad\varphi)\cdot\boldsymbol{\Phi}\d x,\\
\label{VISC_LIM_28c} 0&= \langle\delt\varphi{,}\Phi\rangle_{H^1}  +\intO m(\varphi)\grad\mu\cdot\grad\Phi + (\grad\varphi\cdot\v+\varphi\Gamma_{\v}(\varphi,\sigma)  -\Gamma_{\varphi}(\varphi,\sigma))\Phi\d x,\\
\label{VISC_LIM_28d}0&=\intO (\mu-\epsilon^{-1}\Psi'(\varphi)+ \epsilon\laplace\varphi + \chi\sigma)\Phi\d x,\\
\label{VISC_LIM_28e}0&= \intO \grad\sigma\cdot\grad\Phi + h(\varphi)\sigma\Phi\d x + \int_{\del\Omega}K(\sigma-\sigma_{\infty})\Phi\d\mathcal{H}^{d-1} ,
\end{align}
\end{subequations}
holds for a.\,e. $t\in (0,T)$ and all $\Phi\in H^1,\,\boldsymbol{\Phi}\in\H^1$.
The initial condition is satisfied since $\varphi_{n}(0) = \varphi_0$ a.\,e. in $\Omega$ and by the strong convergence $\varphi_{n}\to\varphi$ in $C^0(L^2)$ as $n\to\infty$.
By the weak (weak-star) lower semi-continuity of norms and \eqref{VISC_LIM_13}, we obtain that $(\varphi,\mu,\sigma,\v,p)$ satisfies
\begin{align}
\nonumber &\norm{\varphi}_{W^{1,\frac{8}{5}}((H^1)^*)\cap L^{\infty}(H^1)\cap L^4(H^2)\cap L^2(H^3)} + \norm{\mu}_{L^2(H^1)} + \norm{\sigma}_{L^2(H^1)} + \norm{\divergence(\varphi\v)}_{L^{\frac{8}{5}}((H^1)^*)}\\
\label{VISC_LIM_29}&\quad + \norm{\divergence(\v)}_{L^2(L^6)} + \norm{\v}_{L^2(\L_{\divergence}^2(\Omega))} + \norm{p}_{L^2(L^2)}\leq C.
\end{align}
\paragraph{Step 5:} Using \eqref{VISC_LIM_28b} and \eqref{VISC_LIM_29}, we obtain that $p$ has a weak derivative in $L^{\frac{8}{5}}(L^2)$ and it holds
\begin{equation}
\label{VISC_LIM_30}\grad p = -\nu\v + (\mu+\chi\sigma)\grad\varphi\quad\text{a.\,e. in }Q.
\end{equation}
By \eqref{LEM_GAGNIR_EST}, we have
\begin{equation*}
\normL{3}{\grad\varphi}\leq C\normL{2}{\grad\varphi}^{\frac{3}{4}}\normH{2}{\grad\varphi}^{\frac{1}{4}}\leq C\normh{1}{\varphi}^{\frac{3}{4}}\normh{3}{\varphi}^{\frac{1}{4}}.
\end{equation*}
Using \eqref{VISC_LIM_29}, this implies
\begin{equation*}
\norm{\grad\varphi}_{L^8(\L^3)}\leq C.
\end{equation*}
Then, using the continuous embedding $H^1\hookrightarrow L^6$ and \eqref{VISC_LIM_29} again, we obtain
\begin{equation*}
\norm{(\mu+\chi\sigma)\grad\varphi}_{L^\frac{8}{5}(\L^2)}\leq C.
\end{equation*}
Since $\v\in L^2(\L^2)$ and $p\in L^2(L^2)$ with bounded norm, from \eqref{VISC_LIM_30} we obtain
\begin{equation}
\label{VISC_LIM_31}\norm{p}_{L^2(L^2)\cap L^\frac{8}{5}(H^1)}\leq C.
\end{equation}
Integrating \eqref{VISC_LIM_28b} by parts, we obtain
\begin{equation*}
-\int_{\del\Omega}p\boldsymbol{\Phi}\cdot\n\d\mathcal{H}^{d-1} = \intO (-\grad p + (\mu+\chi\sigma)\grad\varphi -\nu\v)\cdot\boldsymbol{\Phi}\d x 
\end{equation*}
for all $\boldsymbol{\Phi}\in \H^1$ and a.\,e. $t\in (0,T)$. Because of \eqref{VISC_LIM_30}, this implies
\begin{equation*}
\int_{\del\Omega}p\boldsymbol{\Phi}\cdot\n\d\mathcal{H}^{d-1} = 0
\end{equation*}
for all $\boldsymbol{\Phi}\in \H^1$ and a.\,e. $t\in (0,T)$. Therefore, we obtain
\begin{equation}
\label{VISC_LIM_32}p = 0\quad\text{a.\,e. on }\Sigma.
\end{equation}
With similar arguments, it is straightforward to show that
\begin{equation}
\label{VISC_LIM_33}\mu = \epsilon^{-1}\psi'(\varphi)-\epsilon\laplace\varphi -\chi\sigma\quad\text{a.\,e. in }Q,\qquad \deln\varphi = 0\quad\text{a.\,e. on }\Sigma,
\end{equation}
which completes the proof.
\section{Continuous dependence}
In the following, since it has no bearing on the further analysis, we set $\epsilon = 1$. Let $(\varphi_i,\mu_i,\sigma_i,\v_i,p_i)_{i=1,2}$ be two solutions of \eqref{MEQ}, \eqref{BIC} in the sense of Definition \ref{DEF_WSOL_1}. We denote $\Gamma_{\v}(\varphi_i,\sigma_i)\coloneqq \Gamma_{\v,i},\,\Gamma_{\varphi}(\varphi_i,\sigma_i)\coloneqq \Gamma_{\varphi,i},\,i=1,2,$ and $\sigma_{\infty} \coloneqq \sigma_{1,\infty}-\sigma_{2,\infty}$. Then, the difference $f\coloneqq f_1-f_2$, $f_i\in\{\varphi_i,\mu_i,\sigma_i,\v_i,p_i\},\,i=1,2,$ satisfies

\begin{equation*}
\divergence(\v) = \Gamma_{\v,1}-\Gamma_{\v,2}\quad\text{ a.\,e. in }Q,\quad \varphi(0)=\varphi_{1,0}-\varphi_{2,0}\eqqcolon \varphi_0\quad\text{ a.\,e. in }\Omega,
\end{equation*}
and 
\begin{subequations}
	\label{DIFF}
	\begin{align}
	\nonumber 0 &=\intO (2\eta(\varphi_1)\D\v + \lambda(\varphi_1)\divergence(\v)\I - p\I )\colon \grad\boldsymbol{\Phi}+\nu\v\cdot\boldsymbol{\Phi}\dx  - \intO (\mu+\chi\sigma)\grad\varphi_1\cdot\boldsymbol{\Phi} + (\mu_2+\chi\sigma_2)\grad\varphi\cdot\boldsymbol{\Phi}\dx\\
	\label{DIFF_1}&\quad +\intO (2(\eta(\varphi_1)-\eta(\varphi_2))\D\v_2 +  (\lambda(\varphi_1)-\lambda(\varphi_2))\divergence(\v_2)\I )\colon \grad\boldsymbol{\Phi} \dx, \\
	\nonumber 0 &=\langle\delt\varphi{,}\Phi\rangle_{H^1,(H^1)^*}  +\intO \grad\mu\cdot\grad\Phi + (\varphi_2(\Gamma_{\v,1}-\Gamma_{\v,2})-(\Gamma_{\varphi,1}-\Gamma_{\varphi,2}))\Phi\dx \\
	\label{DIFF_2} &\quad +\intO (\grad\varphi_1\cdot\v + \grad\varphi\cdot\v_2) \Phi  +\varphi\Gamma_{\v,1}\Phi\dx ,\\
	\label{DIFF_4} 0  &= \intO \grad\sigma\cdot\grad\Phi\dx   + \intO  (h(\varphi_1)\sigma + \sigma_2(h(\varphi_1)-h(\varphi_2))\Phi\dx + \int_{\del\Omega}K(\sigma -\sigma_{\infty})\Phi\d\mathcal{H}^{d-1}, 
	\end{align}
\end{subequations}
for a.\,e. $t\in(0,T)$ and for all $\boldsymbol{\Phi}\in \H^1,\,\Phi\in H^1$, where $\mu$ is given by
\begin{equation}
\label{DIFF_3} \mu = \psi'(\varphi_1)-\psi'(\varphi_2)-\laplace\varphi -\chi\sigma.
\end{equation}

\paragraph{Step 1:}  Taking $\Phi = \sigma$ in \eqref{DIFF_4}, we obtain
\begin{equation*}
\intO |\grad\sigma|^2\dx  + K\int_{\del\Omega}|\sigma|^2\d\mathcal{H}^{d-1} + \intO h(\varphi_1)|\sigma|^2\dx  = -\intO \sigma_2(h(\varphi_1)-(\varphi_2))\sigma\dx  + K\int_{\del\Omega}\sigma\sigma_{\infty}\d\mathcal{H}^{d-1}.
\end{equation*}
Using the non-negativity of $h(\cdot)$, we can neglect the third term on the l.h.s. of this equation. Applying Hölder's and Young's inequalities, we therefore obtain
\begin{equation}
\label{UNIQUE_1}\intO |\grad\sigma|^2\dx  + K\int_{\del\Omega}|\sigma|^2\d\mathcal{H}^{d-1} \leq \frac{L_h^2}{4\delta}\norml{3}{\sigma_2}^2\norml{2}{\varphi}^2 + \frac{K}{2}\left(\norm{\sigma}_{L^2(\del\Omega)}^2 + \norm{\sigma_{\infty}}_{L^2(\del\Omega)}^2\right) +\delta\norml{6}{\sigma}^2,
\end{equation}
with $\delta >0$ to be chosen and where we used (B2). Using Poincaré's inequality and the continuous embedding $H^1\hookrightarrow L^6$, we have
\begin{equation}
\label{UNIQUE_2}\norml{6}{\sigma}^2\leq 2C_P^2\left(\normL{2}{\grad\sigma}^2 + \norm{\sigma}_{L^2(\del\Omega)}^2\right).
\end{equation}
Choosing 
\begin{equation*}
\delta = \frac{1}{8C_P^2}\min\{1,K\},
\end{equation*}
and using \eqref{UNIQUE_2}, this implies 
\begin{equation}
\label{UNIQUE_3}\normh{1}{\sigma} \leq C(K,L_h,\Omega)\left(\norml{6}{\sigma_2}\norml{2}{\varphi}+\norm{\sigma_{\infty}}_{L^2(\del\Omega)}\right).
\end{equation}
\paragraph{Step 2:} With similar arguments as in Section 3, we deduce the existence of a solution $\u \in \H^1$ of the problem
\begin{equation*}
\divergence(\u ) = \Gamma_{\v,1}-\Gamma_{\v,2} \quad\text{in }\Omega,\quad
\u  = \frac{1}{|\del\Omega|}\left(\intO\Gamma_{\v,1}-\Gamma_{\v,2}\dx \right)\n \quad\text{on }\del\Omega,
\end{equation*}
satisfying the estimate
\begin{equation}
\label{UNIQUE_4}\normH{1}{\u }\leq c\norml{2}{\Gamma_{\v,1}-\Gamma_{\v,2}},
\end{equation}
with a constant $c$ depending only on $\Omega$. Choosing $\boldsymbol{\Phi} = \v-\u $ in \eqref{DIFF_1} and $\Phi = \varphi -\laplace\varphi$ in \eqref{DIFF_2}, integrating by parts, using  \eqref{DIFF_3} and summing the resulting identities, we obtain
\begin{align}
\nonumber &\frac{\d}{\dt }\frac{1}{2}(\norm{\grad\varphi}_{\L^2}^2 + \norml{2}{\varphi}^2) + \intO 2\eta|\D\v|^2 + \nu|\v|^2 + |\laplace\varphi|^2+ |\grad\laplace\varphi|^2\dx  \\
\nonumber&= \intO \grad(\psi'(\varphi_1)-\psi'(\varphi_2))\cdot(\grad\laplace\varphi -\grad\varphi) - ((\Gamma_{\varphi,1}-\Gamma_{\varphi,2})-\varphi_2(\Gamma_{\v,1}-\Gamma_{\v,2}))\laplace\varphi\dx \\
\nonumber&\quad +\intO (\grad\varphi\cdot\v_2  +\varphi\Gamma_{\v,1})\laplace\varphi + (\psi'(\varphi_1)-\psi'(\varphi_2))\grad\varphi_1\cdot(\v-\u )\dx,\\
\nonumber&\quad + \intO  (\mu_2+\chi\sigma_2)\grad\varphi\cdot(\v-\u ) + \laplace\varphi\grad\varphi_1\cdot\u \dx +\intO \chi\grad\sigma\cdot(\grad\varphi-\grad\laplace\varphi)\d x \\
\nonumber&\quad + \intO ((\Gamma_{\varphi,1}-\Gamma_{\varphi,2})-\varphi_2(\Gamma_{\v,1}-\Gamma_{\v,2}))\varphi -(\grad\varphi_1\cdot\v + \grad\varphi\cdot\v_2) \varphi  -\Gamma_{\v,1}|\varphi|^2\dx \\
\label{UNIQUE_5}&\quad + \intO 2\eta(\varphi_1) \D\v\colon \grad\u  + \nu\v\cdot\u  -2(\eta(\varphi_1)-\eta(\varphi_2))\D\v_2\colon \grad(\v-\u )\dx.
\end{align}
\paragraph{Step 3:} We now estimate the terms on the r.h.s. of \eqref{UNIQUE_5}. In the following, we will frequently use Young's, Hölder's inequalities and Lemma \ref{LEM_GAGNIR}. First of all, using \textnormal{(A4)} and (B2), we obtain
\begin{equation*}
|\Gamma_{\v,1}-\Gamma_{\v,2}|\leq   C(|\sigma| + |\sigma_2||\varphi| + |\varphi|),\quad
|\Gamma_{\varphi,1}-\Gamma_{\varphi,2}|\leq  C(|\sigma| + |\sigma_2||\varphi| + |\varphi|).
\end{equation*}
Hence, applying \eqref{UNIQUE_3} and the continuous embedding $H^1\hookrightarrow L^6$ yields
\begin{equation}
\label{UNIQUE_6}\norml{2}{\Gamma_{\v,1}-\Gamma_{\v,2}}+ \norml{2}{\Gamma_{\varphi,1}-\Gamma_{\varphi,2}}\leq C\left(1+\norml{6}{\sigma_2}\right) \normh{1}{\varphi}+C\norm{\sigma_{\infty}}_{L^2(\del\Omega)}.
\end{equation}
Using \eqref{THM_WSOL_1_EST_1}, \eqref{UNIQUE_6} and the continuous embedding $H^1\hookrightarrow L^6$, we obtain
\begin{align}
\nonumber&\left|\intO ((\Gamma_{\varphi,1}-\Gamma_{\varphi,2})-\varphi_2(\Gamma_{\v,1}-\Gamma_{\v,2}))\laplace\varphi\dx \right|\\
\nonumber &\quad\leq C\big(1+\norml{\infty}{\varphi_2}^2\big)\big((1+\norml{6}{\sigma_2}^2)\normh{1}{\varphi}^2 +\norm{\sigma_{\infty}}_{L^2(\del\Omega)}^2\big) + \frac{1}{8}\norml{2}{\laplace\varphi}^2\\
\label{UNIQUE_7} &\quad\leq C\left(1+\normh{2}{\varphi_2}+\normh{1}{\sigma_2}^4\right)\normh{1}{\varphi}^2+ C\big(1+\normh{2}{\varphi_2}\big)\norm{\sigma_{\infty}}_{L^2(\del\Omega)}^2 + \frac{1}{8}\norml{2}{\laplace\varphi}^2 .
\end{align}
By the specific form of $\Gamma_{\v}$ and \textnormal{(A4)}, applying the continuous embedding $H^1\hookrightarrow L^6$ gives
\begin{equation}
\label{UNIQUE_8}\left|\intO \varphi\Gamma_{\v,1}\laplace\varphi\dx \right|\leq C\left(1+\norml{3}{\sigma_1}^2\right)\normh{1}{\varphi}^2 + \frac{1}{8}\norml{2}{\laplace\varphi}^2.
\end{equation}
Hence, recalling the continuous embedding $\H^1\hookrightarrow \L^6$ and \eqref{ELLIPTIC_EST}, we calculate
\begin{align}
\nonumber \left|\intO \grad\varphi\cdot\v_2\laplace\varphi\dx \right| &\leq \normL{3}{\grad\varphi}\normL{6}{\v_2}\norml{2}{\laplace\varphi}\\
\nonumber &\leq C\normH{1}{\v_2}\normL{2}{\grad\varphi}\normh{2}{\varphi}^\frac{1}{2}\normL{2}{\grad\laplace\varphi}^{\frac{1}{2}}\\
\nonumber &\leq C\normH{1}{\v_2}\normh{1}{\varphi}\left(\norml{2}{\varphi} + \norml{2}{\laplace\varphi}+\normL{2}{\grad\laplace\varphi}\right)\\
\label{UNIQUE_9}&\leq C\left(1+\normH{1}{\v_2}^2\right)\normh{1}{\varphi}^2 + \frac{1}{8}\norml{2}{\laplace\varphi}^2 + \frac{1}{4}\normL{2}{\grad\laplace\varphi}^2.
\end{align}
Due to \textnormal{(B3)}, \eqref{LEM_GAGNIR_EST}, \eqref{THM_WSOL_1_EST_1} and the continuous embedding $H^1\hookrightarrow L^6$, we obtain
\begin{equation}
\label{UNIQUE_10} \norml{2}{\psi'(\varphi_1)-\psi'(\varphi_2)}^2\leq  C\left(1+\normh{2}{\varphi_1}^2+\normh{2}{\varphi_2}^2\right)\normh{1}{\varphi}^2.
\end{equation}
By the continuous embedding $\H^1\hookrightarrow \L^6$ and \eqref{THM_WSOL_1_EST_1}, this gives
\begin{align}
\nonumber \left|\intO (\psi'(\varphi_2)-\psi'(\varphi_1))\grad\varphi_1\cdot\v\dx \right|&\leq C\norml{2}{\psi'(\varphi_1)-\psi'(\varphi_2)}^2\normL{3}{\grad\varphi_1}^2 + \delta_1\normL{6}{\v}^2\\
\label{UNIQUE_11}&\leq C\left(1+\normh{2}{\varphi_1}^4+\normh{2}{\varphi_2}^4\right)\normh{1}{\varphi}^2 + C\delta_1\normH{1}{\v}^2,
\end{align}
with $\delta_1>0$ to be chosen later.
Recalling \eqref{UNIQUE_4} and \eqref{UNIQUE_6}, we obtain with similar arguments that
\begin{align}
\nonumber \left|\intO (\psi'(\varphi_2)-\psi'(\varphi_1))\grad\varphi_1\cdot\u \dx \right|&\leq C\norml{2}{\psi'(\varphi_1)-\psi'(\varphi_2)}\normL{3}{\grad\varphi_1}\normL{6}{\u}\\
\nonumber &\leq C\left(1+\norml{12}{\varphi_1}^4+\norml{12}{\varphi_2}^4\right)\normh{1}{\varphi}\normL{3}{\grad\varphi_1}\normH{1}{\u}\\
\label{UNIQUE_12} &\leq C\left(1+\normh{2}{\varphi_1}^4+\normh{2}{\varphi_2}^4 + \norml{6}{\sigma_2}^4\right)\normh{2}{\varphi}^2  + C\norm{\sigma_{\infty}}_{L^2(\del\Omega)}^2 .
\end{align}
Again using \eqref{UNIQUE_4} and \eqref{UNIQUE_6} and the continuous embeddings $H^1\hookrightarrow L^6,\,\H^1\hookrightarrow \L^6$, we obtain
\begin{equation}
\label{UNIQUE_13} \left|\intO (\mu_2+\chi\sigma_2)\grad\varphi\cdot(\v-\u )\dx \right|\leq C\left(1+\normh{1}{\mu_2}^2+\normh{1}{\sigma_2}^2\right)\normh{1}{\varphi}^2 + C\norm{\sigma_{\infty}}_{L^2(\del\Omega)}^2 + \delta_2\normH{1}{\v}^2,
\end{equation}
with $\delta_2>0$ to be chosen later. Once more using \eqref{THM_WSOL_1_EST_1},\eqref{UNIQUE_4} and \eqref{UNIQUE_6} and the continuous embedding $\H^1\hookrightarrow \L^6$, we deduce that
\begin{align}
\nonumber\left|\intO \laplace\varphi\grad\varphi_1\cdot\u \dx \right|&\leq \norml{2}{\laplace\varphi}\normL{3}{\grad\varphi_1}\normL{6}{\u}\\
&\nonumber \leq C\normL{3}{\grad\varphi_1}^2\normH{1}{\u}^2 + \frac{1}{8}\norml{2}{\laplace\varphi}^2\\
\label{UNIQUE_14}&\leq  C\left(1+\normh{2}{\varphi_1}^2 + \normh{1}{\sigma_2}^4\right)\normh{1}{\varphi}^2 + C\normh{2}{\varphi_1}\norm{\sigma_{\infty}}_{L^2(\del\Omega)}^2+ \frac{1}{8}\norml{2}{\laplace\varphi}^2.
\end{align}
Using the assumptions on $\eta(\cdot)$, we obtain
\begin{align}
\nonumber\left|\intO 2\eta(\varphi_1) \D\v\colon\grad\u  + \nu\v\cdot\u \dx \right|&\leq \frac{\eta_0}{2}\normL{2}{\D\v}^2 + \frac{\nu}{2}\normL{2}{\v}^2 + C\normH{1}{\u}^2\\
\nonumber&\leq \frac{\eta_0}{2}\normL{2}{\D\v}^2 + \frac{\nu}{2}\normL{2}{\v}^2 + C\left(1+\normh{1}{\sigma_2}^2\right)\normh{1}{\varphi}^2\\
\label{UNIQUE_15}&\quad +C\norm{\sigma_{\infty}}_{L^2(\del\Omega)}^2 .
\end{align}
Furthermore, using the boundedness of $\eta'(\cdot)$, elliptic regularity and \eqref{LEM_GAGNIR_EST}, \eqref{UNIQUE_4}, \eqref{UNIQUE_6}, we obtain
\begin{align}
\nonumber &\left|\intO 2(\eta(\varphi_1)-\eta(\varphi_2)\D\v_2\colon\grad(\v-\u )\dx \right| \\
\nonumber &\leq C\norml{\infty}{\varphi}\normL{2}{\D\v_2}\left(\normL{2}{\grad\v}+\normL{2}{\grad\u }\right)\\
\nonumber &\leq C\normh{1}{\varphi}^{\frac{3}{4}}\normh{3}{\varphi}^{\frac{1}{4}}\normL{2}{\D\v_2}\left(\normL{2}{\grad\v}+\normL{2}{\grad\u}\right)\\
\nonumber &\leq \delta_3\normL{2}{\grad\v}^2+C\normL{2}{\grad\u}^2 + C\normh{1}{\varphi}^{\frac{3}{2}}\normh{3}{\varphi}^{\frac{1}{2}}\normL{2}{\D\v_2}^2\\
&\leq \delta_3\normL{2}{\grad\v}^2+C\norm{\sigma_{\infty}}_{L^2(\del\Omega)}^2+ \frac{1}{4}\normL{2}{\grad\laplace\varphi}^2 + C\left(1+\normL{2}{\D\v_2}^{\frac{8}{3}}+\normh{1}{\sigma_2}^2\right)\normh{1}{\varphi}^2,
\end{align}
with $\delta_3>0$ to be chosen later.
Due to \eqref{UNIQUE_6} and since $\varphi_2\in L^{\infty}(H^1)$ with bounded norm, we get
\begin{equation}
\label{UNIQUE_16}\left|\intO (1-\varphi_2)(\Gamma(\varphi_1,\sigma_1)-\Gamma(\varphi_2,\sigma_2))\varphi\dx \right|\leq C\left(1+\norm{\sigma_2}_{L^6}\right)\normh{1}{\varphi}^2+C\norm{\sigma_{\infty}}_{L^2(\del\Omega)}^2.
\end{equation}
Applying \eqref{THM_WSOL_1_EST_1} and the continuous embeddings $H^1\hookrightarrow L^6,\,\H^1\hookrightarrow\L^6$ gives
\begin{align}
\nonumber \left|\intO (\grad\varphi_1\cdot\v + \grad\varphi\cdot\v_2)\varphi\dx \right| &\leq \norm{\grad\varphi_1}_{\L^3}\norm{\v}_{\L^6}\norml{2}{\varphi} + \norm{\grad\varphi}_{\L^2}\norm{\v_2}_{\L^6}\norm{\varphi}_{L^3}\\
\nonumber &\leq C\left(\norm{\grad\varphi_1}_{\L^3}\norm{\v}_{\H^1}\normh{1}{\varphi} + \normh{1}{\varphi}^2\norm{\v_2}_{\H^1}\right)\\
\label{UNIQUE_17}&\leq C\left(\normh{2}{\varphi_1} + \norm{\v_2}_{\H^1}\right)\normh{1}{\varphi}^2 + \delta_4\norm{\v}_{\H^1}^2,
\end{align}
with $\delta_4>0$ to be chosen later. For the last term on the r.h.s. of \eqref{UNIQUE_5}, we use \textnormal{(A4)} and the continuous embedding $H^1\hookrightarrow L^6$ to obtain
\begin{equation}
\label{UNIQUE_18}\left|\intO \Gamma_{\v,1}|\varphi|^2\dx \right|\leq C\left(1+\norml{6}{\sigma_2}\right)\normh{1}{\varphi}^2.
\end{equation}
We now estimate the first term on the r.h.s. of \eqref{UNIQUE_5}.
To this end, we first observe that
\begin{equation*}
\grad(\psi'(\varphi_1)-\psi'(\varphi_2)) = \psi''(\varphi_1)\grad\varphi_1 - \psi''(\varphi_2)\grad\varphi_2 = \psi''(\varphi_1)\grad\varphi + \grad\varphi_2(\psi''(\varphi_1)-\psi''(\varphi_2)).
\end{equation*}
Due to \eqref{ASS_PSI_2}-\eqref{ASS_PSI_3} and \eqref{THM_WSOL_1_EST_1}, we obtain
\begin{equation*}
\intO |\psi''(\varphi_1)\grad\varphi|^2\dx  \leq C\intO (1+|\varphi_1|^8)|\grad\varphi|^2\dx \leq C\left(1+\norml{\infty}{\varphi_1}^8\right)\normL{2}{\grad\varphi}^2\leq C\left(1+\normh{3}{\varphi_1}^2\right)\normh{1}{\varphi}^2.
\end{equation*}
Applying \eqref{THM_WSOL_1_EST_1} and \eqref{ASS_PSI_7}, we conclude that
\begin{align*}
\intO |\grad\varphi_2(\psi''(\varphi_1)-\psi''(\varphi_2))|^2\dx &\leq C\intO (1+|\varphi_1|^6+|\varphi_2|^6)|\grad\varphi_2|^2|\varphi|^2\dx \\
&\leq C\left(1+\norml{18}{\varphi_1}^6 +\norml{18}{\varphi_2}^6\right)\normL{6}{\grad\varphi_2}^2\norml{6}{\varphi}^2\\
&\leq C\left(1+\normh{2}{\varphi_1}^4 +\normh{2}{\varphi_2}^4\right)\normh{1}{\varphi}^2.
\end{align*}
Combining the last two estimates, we obtain
\begin{equation}
\label{UNIQUE_19}\norm{\grad(\psi'(\varphi_1)-\psi'(\varphi_2))}_{\L^2}^2\leq C\left(1+\normh{2}{\varphi_1}^4 +\normh{2}{\varphi_2}^4 + \normh{3}{\varphi_1}^2\right)\normh{1}{\varphi}^2.
\end{equation}
From this inequality, we deduce that
\begin{align}
\nonumber \left|\intO \grad(\psi'(\varphi_1)-\psi'(\varphi_2))\cdot(\grad\laplace\varphi-\grad\varphi)\dx \right|&\leq C\left(1+\normh{2}{\varphi_1}^4 +\normh{2}{\varphi_2}^4 + \normh{3}{\varphi_1}^2\right)\normh{1}{\varphi}^2\\
\label{UNIQUE_20}&\quad  + \frac{1}{8}\normL{2}{\grad\laplace\varphi}^2.
\end{align}
Finally, using \eqref{UNIQUE_3}, it is easy to check that
\begin{equation}
\label{UNIQUE_20a}\left|\intO \chi\grad\sigma\cdot(\grad\varphi-\grad\laplace\varphi)\d x\right|\leq C\left(1+\norml{6}{\sigma_2}^2\right)\normh{1}{\varphi}^2 + C\norm{\sigma_{\infty}}_{L^2(\del\Omega)}^2 + \frac{1}{8}\normL{2}{\grad\laplace\varphi}^2.
\end{equation}
Using \eqref{UNIQUE_7}-\eqref{UNIQUE_20a} in \eqref{UNIQUE_5} and choosing
\begin{equation*}
\delta_1 = \delta_2 =\delta_3= \delta_4= \frac{\min\{\frac{\nu}{2},\eta\}}{8C_K^2}\eqqcolon C_1,
\end{equation*}
where $C_K$ is Korn's constant, we end up with
\begin{equation}
\label{UNIQUE_21}\frac{\d}{\dt }\frac{1}{2}\left(\normL{2}{\grad\varphi}^2 + \norml{2}{\varphi}^2\right) + 4C_1\normH{1}{\v}^2 + \frac{1}{2}\intO |\laplace\varphi|^2+ |\grad\laplace\varphi|^2\dx  \leq \alpha_1(t)\normh{1}{\varphi}^2+ \alpha_2(t)\norm{\sigma_{\infty}}_{L^2(\del\Omega)}^2,
\end{equation}
where 
\begin{align*}
\alpha_1(t)&\coloneqq C\left(1+\normh{2}{\varphi_1}^4 + \normh{2}{\varphi_2}^4 + \normh{3}{\varphi_1}^2+\normh{2}{\varphi_2}^2+\normH{1}{\v_2}^{\frac{8}{3}} + \normh{1}{\mu_2}^2+\normh{1}{\sigma_1}^2+\normh{1}{\sigma_2}^4\right)\\
\alpha_2(t)&\coloneqq C\left(1+\normh{2}{\varphi_1} + \normh{2}{\varphi_2}\right).
\end{align*}
Due to \eqref{THM_WSOL_1_EST_1}-\eqref{THM_WSOL_1_EST_2}, it follows that $\alpha_1\in L^1(0,T),\,\alpha_2\in L^4(0,T)$. We remark that $\alpha_1\in L^1(0,T)$ only holds provided $\sigma_2\in L^4(H^1)$ with bounded norm. Then, using a Gronwall argument (see \cite[Lemma 3.1]{GarckeLam3}) in \eqref{UNIQUE_21} yields
\begin{equation*}
\sup_{s\in (0,T]}\normh{1}{\varphi(s)}^2 + \int_{0}^{T}\normH{1}{\v}^2\d s + \int_{0}^{T}\intO |\laplace\varphi|^2+ |\grad\laplace\varphi|^2\dx \d s \leq C\left(\normh{1}{\varphi_0}^2 + \norm{\sigma_{\infty}}_{L^4(L^2(\del\Omega))}^2\right). 
\end{equation*}
Together with elliptic regularity theory, this gives
\begin{equation}
\label{UNIQUE_22}\norm{\varphi}_{L^{\infty}(H^1)\cap L^2(H^3)} + \norm{\v}_{L^2(\H^1)} \leq C\left( \normh{1}{\varphi_0} + \norm{\sigma_{\infty}}_{L^4(L^2(\del\Omega))}\right).
\end{equation}
Now, from \eqref{UNIQUE_3} and \eqref{UNIQUE_22}, we immediately obtain
\begin{equation}
\label{UNIQUE_23}\norm{\sigma}_{L^2(H^1)}\leq C\left(\normh{1}{\varphi_0} + \norm{\sigma_{\infty}}_{L^4(L^2(\del\Omega))}\right).
\end{equation}
Using \eqref{LEM_GAGNIR_EST}, \eqref{UNIQUE_10}, \eqref{UNIQUE_19}, \eqref{UNIQUE_22} and the boundedness of $\varphi_1,\varphi_2\in L^{\infty}(H^1)\cap L^4(H^2)\cap L^2(H^3)$, it is straightforward to check that
\begin{equation*}
\norm{\psi'(\varphi_1)-\psi'(\varphi_2)}_{L^2(H^1)}\leq C\left(\normh{1}{\varphi_0} + \norm{\sigma_{\infty}}_{L^4(L^2(\del\Omega))}\right).
\end{equation*}
Together with \eqref{UNIQUE_22}-\eqref{UNIQUE_23}, a comparison argument in \eqref{DIFF_3} yields
\begin{equation}
\label{UNIQUE_24} \norm{\mu}_{L^2(H^1)}  \leq C\left( \normh{1}{\varphi_0} + \norm{\sigma_{\infty}}_{L^4(L^2(\del\Omega))}\right).
\end{equation}
Using \eqref{UNIQUE_22}-\eqref{UNIQUE_24}, a comparison argument in \eqref{DIFF_2} yields
\begin{equation}
\label{UNIQUE_25} \norm{\varphi}_{H^1((H^1)^*)}  \leq C\left( \normh{1}{\varphi_0} + \norm{\sigma_{\infty}}_{L^4(L^2(\del\Omega))}\right).
\end{equation}
Combining \eqref{UNIQUE_22}-\eqref{UNIQUE_25}, we deduce the estimate
\begin{align}
\nonumber&\norm{\varphi}_{H^1((H^1)^*)\cap L^{\infty}(H^1)\cap L^2(H^3)} + \norm{\mu}_{L^2(H^1)} + \norm{\sigma}_{L^2(H^1)} + \norm{\v}_{L^2(\H^1)}\\
\label{UNIQUE_26}&\quad \leq C\left( \normh{1}{\varphi_0} + \norm{\sigma_{\infty}}_{L^4(L^2(\del\Omega))}\right).
\end{align}
\paragraph{Step 4:} It remains to get an estimate on the pressure. To this end, let $\q\in \H^1$ be a solution of
\begin{equation*}
\divergence(\q) = p\quad\text{in } \Omega,\quad \q = \frac{1}{|\del\Omega|}\left(\intO p\dx \right)\n \quad\text{on }\del\Omega,
\end{equation*}
such that
\begin{equation}
\label{UNIQUE_27}\normH{1}{\q}\leq c\norml{2}{p},
\end{equation}
with $c$ depending only on $\Omega$. Then, choosing $\boldsymbol{\Phi} = \q$ in \eqref{DIFF_1}, we obtain
\begin{align}
\nonumber \norml{2}{p}^2 &= \intO (2\eta(\varphi_1)\D\v + \lambda(\varphi_1)\divergence(\v)\I)\colon \grad\q+\nu\v\cdot\q\dx  - \intO (\mu+\chi\sigma)\grad\varphi_1\cdot\q + (\mu_2+\chi\sigma_2)\grad\varphi\cdot\q\dx\\
\label{UNIQUE_28}&\quad +\intO (2(\eta(\varphi_1)-\eta(\varphi_2))\D\v_2 +  (\lambda(\varphi_1)-\lambda(\varphi_2))\divergence(\v_2)\I )\colon \grad\q \dx.
\end{align}
Using \eqref{THM_WSOL_1_EST_1}-\eqref{THM_WSOL_1_EST_2} and (A3), a straightforward calculation shows that
\begin{align*}
&\left|\intO (2\eta(\varphi_1)\D\v + \lambda(\varphi_1)\divergence(\v)\I)\colon \grad\q+\nu\v\cdot\q\dx  - \intO (\mu+\chi\sigma)\grad\varphi_1\cdot\q + (\mu_2+\chi\sigma_2)\grad\varphi\cdot\q\dx\right|\\
&\leq C\left(\normH{1}{\v}^2 + \normh{1}{\mu+\chi\sigma}^2 + \normh{1}{\mu_2+\chi\sigma_2}^2\normh{1}{\varphi}^2\right) + \frac{1}{4}\norml{2}{p}^2.
\end{align*}
For the remaining terms, we use again \eqref{THM_WSOL_1_EST_1}-\eqref{THM_WSOL_1_EST_2} and (A3) to obtain
\begin{equation*}
\left|\intO (2(\eta(\varphi_1)-\eta(\varphi_2))\D\v_2 +  (\lambda(\varphi_1)-\lambda(\varphi_2))\divergence(\v_2)\I )\colon \grad\q \dx\right|\\
\leq C\normH{1}{\v_2}^2\norml{\infty}{\varphi}^2  + \frac{1}{4}\norml{2}{p}^2.
\end{equation*}
Using the last two inequalities in \eqref{UNIQUE_28}, integrating the resulting estimate in time from $0$ to $T$ and using Young's generalised inequality, we deduce that
\begin{equation*}
\norm{p}_{L^2(L^2)}^2\leq C\left(\norm{\v}_{L^2(\H^1)}^2 + \norm{\mu+\chi\sigma}_{L^2(H^1)}^2 + \norm{\mu_2+\chi\sigma_2}_{L^2(H^1)}^2\norm{\varphi}_{L^{\infty}(H^1)}^2 + \norm{\v_2}_{L^{\frac{8}{3}}(\H^1)}^2\norm{\varphi}_{L^8(L^{\infty})}^2\right).
\end{equation*}
Therefore, using \eqref{THM_WSOL_1_EST_1}-\eqref{THM_WSOL_1_EST_2}, \eqref{UNIQUE_26} and the continuous embedding $L^{8}(L^{\infty})\hookrightarrow L^{\infty}(H^1)\cap L^2(H^3)$, the last inequality implies
\begin{equation}
\label{UNIQUE_29}\norm{p}_{L^2(L^2)}\leq C\left(\normh{1}{\varphi_0} + \norm{\sigma_{\infty}}_{L^4(L^2(\del\Omega))}\right).
\end{equation}
Together with \eqref{UNIQUE_26}, we obtain \eqref{THM_CONTDEP_EST}, whence the proof is complete.  
\section{Existence of strong solutions}
In this section, we will prove Theorem \ref{THM_SSOL}. The testing procedure can again be justified by a Galerkin scheme. 
In the following, we assume for simplicity and as it has no further consequence for the analysis that $\epsilon = 1$. Then, with similar arguments as before, we can arrive at
\begin{align}
\nonumber&\norm{\varphi}_{H^1((H^1)^*)\cap L^{\infty}(H^1)\cap L^4(H^2)\cap L^2(H^3) } + \norm{\sigma}_{L^4(H^1)} + \norm{\mu}_{L^2(H^1)\cap L^4(L^2)}\\
\label{APRI_SSOL_1}&\quad + \norm{\divergence(\varphi\v)}_{L^2(L^2)}  + \norm{\v}_{L^{\frac{8}{3}}(\H^1)}+ \norm{p}_{L^2(L^2)}\leq C.
\end{align}
We will now show the result in a series of higher order estimates.\newline
\paragraph{Step 1:}
Observing that \eqref{WFORM_1d} is the weak formulation of
\begin{alignat*}{3}
-\laplace\sigma + h(\varphi)\sigma &= 0&&\quad\text{a.~e. in }\Omega,\\
\deln\sigma + K\sigma &= K\sigma_{\infty}&& \quad \text{a.~e. on }\del\Omega,
\end{alignat*}
due to the assumptions on $h(\cdot)$ and using \cite[Thm. 2.4.2.6]{Grisvard} we deduce that
\begin{equation*}
\normh{2}{\sigma}\leq C\norm{K\sigma_{\infty}}_{H^{\frac{1}{2}}(\del\Omega)}.
\end{equation*}
Therefore, using the continuous embedding $H^2\hookrightarrow L^{\infty}$ and the fact that $\sigma_{\infty}\in H^1(H^{\frac{1}{2}}(\del\Omega))\hookrightarrow C^0(H^{\frac{1}{2}}(\del\Omega))$, we have
\begin{equation}
\label{APRI_SSOL_3}\norm{\sigma}_{L^{\infty}(H^2)\cap L^{\infty}(L^{\infty})}\leq C.
\end{equation}
By \textnormal{(A3)}, this yields
\begin{equation}
\label{APRI_SSOL_4}\norm{\divergence(\v)}_{L^{\infty}(L^{\infty})} +  \norm{\Gamma_{\varphi}}_{L^{\infty}(L^{\infty})}\leq C.
\end{equation}
Now, for $h>0$ we introduce the incremental ratio
\begin{equation*}
\del_t^h u(t) = \frac{1}{h}(u(t+h)-u(t)).
\end{equation*}
Then, using \eqref{WFORM_1d}, we see that
\begin{equation*}
0 = \intO \nabla \del_t^h\sigma(t)\cdot\grad\Phi + \big(\del_t^h(h(\varphi(t)))\sigma(t+h) + \del_t^h\sigma(t)h(\varphi(t)))\Phi\dx + \int_{\del\Omega}K(\del_t^h\sigma(t)-\del_t^h\sigma_{\infty}(t))\Phi\d\mathcal{H}^{d-1}
\end{equation*}
holding for almost every $t\in (0,T-h]$. Testing this equation with $\Phi = \del_t^h\sigma(t)$, integrating in time from $0$ to $T-h$ and using the non-negativity of $h(\cdot)$, we obtain
\begin{align}
\nonumber\int_{0}^{T-h}\normL{2}{\grad\del_t^h\sigma(t)}^2\d t + K\int_{0}^{T-h}\norm{\del_t^h\sigma(t)}_{L^2(\del\Omega)}^2\d t&\leq \int_{0}^{T-h}\intO \del_t^h(h(\varphi(t)))\sigma(t+h)\del_t^h\sigma(t)\d x\d t\\
\label{APRI_SSOL_5} &\quad + \int_{0}^{T-h}\int_{\del\Omega}K\del_t^h\sigma(t)\del_t^h\sigma_{\infty}(t)\d\mathcal{H}^{d-1}\d t.
\end{align}
To estimate the r.h.s. of this equation, we use \eqref{APRI_SSOL_3} and the Lipschitz-continuity of $h(\cdot)$ to obtain
\begin{align}
\nonumber\left|\int_{0}^{T-h}\intO \del_t^h(h(\varphi(t)))\sigma(t+h)\del_t^h\sigma(t)\d x\d t\right| &\leq C\int_{0}^{T-h}\intO |\del_t^h\varphi(t)\del_t^h\sigma(t)|\d x\d t\\
\nonumber&\leq C\norm{\del_t^h\varphi}_{L^2(0,T-h;(H^1)^*)}\norm{\del_t^h\sigma}_{L^2(0,T-h;H^1)}\\
\label{APRI_SSOL_6}&\leq C\norm{\del_t\varphi}_{L^2(0,T;(H^1)^*)}\norm{\del_t^h\sigma}_{L^2(0,T-h;H^1)}.
\end{align}
With similar arguments and using the trace theorem, the remaining term can be estimated by
\begin{equation}
\label{APRI_SSOL_6a}\left|\int_{0}^{T-h}\int_{\del\Omega}K\del_t^h\sigma(t)\del_t^h\sigma_{\infty}(t)\d\mathcal{H}^{d-1}\d t\right|\leq C\norm{\del_t\sigma_{\infty}}_{L^2(0,T;H^{\frac{1}{2}}(\del\Omega))}\norm{\del_t^h\sigma}_{L^2(0,T-h;H^1)}.
\end{equation}
Using \eqref{APRI_SSOL_5}-\eqref{APRI_SSOL_6a} together with \eqref{APRI_SSOL_1} and \textnormal{(C2)}, an application of Poincaré's inequality implies
\begin{equation}
\label{APRI_SSOL_7}\norm{\del_t^h\sigma}_{L^2(0,T-h;H^1)}\leq C\left(\norm{\del_t\varphi}_{L^2(0,T;(H^1)^*)}+\norm{\del_t\sigma_{\infty}}_{L^2(0,T;H^{\frac{1}{2}})}\right)\leq C.
\end{equation}
Since the constant $C$ is independent of $h>0$, this yields
\begin{equation*}
\norm{\del_t\sigma}_{L^2(H^1)}\leq C.
\end{equation*}
Combining this inequality with \eqref{APRI_SSOL_3} and using the continuous embedding $H^1(H^1)\hookrightarrow C^0(H^1)$, we obtain
\begin{equation}
\label{APRI_SSOL_8}\norm{\sigma}_{H^1(H^1)\cap C^0(H^1)\cap L^{\infty}(H^2)}\leq C.
\end{equation}
\paragraph{Step 2:}
Choosing  $\Phi = \del_t\varphi$ in \eqref{WFORM_1b} and $\Phi = \laplace\del_t\varphi$ in \eqref{WFORM_1c} and integrating by parts, we obtain
\begin{align}
\nonumber \frac{\d}{\dt }\frac{1}{2}\intO |\laplace\varphi|^2\dx  + \intO |\del_t\varphi|^2 &= -\intO (\divergence(\varphi\v) -\Gamma_{\varphi})\del_t\varphi \dx + \chi\intO \nabla\sigma\cdot\nabla\del_t\varphi\d x  \\
\label{APRI_SSOL_9}&\quad  + \intO \psi'''(\varphi)|\grad\varphi|^2\del_t\varphi + \psi''(\varphi)\laplace\varphi\del_t\varphi\dx .
\end{align}
We recall that $\Gamma_{\varphi},~\Gamma_{\v}\in L^2(L^2)$ with bounded norm. Then, using Hölder's and Young's inequalities, we can estimate the first three terms on the r.h.s. of  \eqref{APRI_SSOL_9} by
\begin{equation}
\label{APRI_SSOL_10}\left|\intO (\divergence(\varphi\v)  -\Gamma_{\varphi} )\delt\varphi \dx \right|\leq  C\left(\norml{2}{\divergence(\varphi\v)}^2 +\norml{2}{\Gamma_{\varphi}}^2\right) +\frac{1}{4}\norml{2}{\delt\varphi}^2.
\end{equation}
For the last term on the r.h.s. of \eqref{APRI_SSOL_9}, we use Hölder's and Young's inequalities together with \eqref{LEM_GAGNIR_EST}, \eqref{ASS_PSI_2}-\eqref{ASS_PSI_3} and \eqref{APRI_SSOL_1} to obtain
\begin{equation}
\label{APRI_SSOL_11} \left|\intO  \psi''(\varphi)\laplace\varphi\delt\varphi\dx \right|\leq  C\left(1+\norml{\infty}{\varphi}^4\right)\norml{2}{\laplace\varphi}\norml{2}{\delt\varphi}\leq \frac{1}{4}\norml{2}{\delt\varphi}^2 + C\left(1+\normh{3}{\varphi}^2\right)\norml{2}{\laplace\varphi}^2.
\end{equation}
Now, using \textnormal{(C3)}, Hölder's and Young's inequalities, \eqref{ELLIPTIC_EST}, \eqref{LEM_GAGNIR_EST} and \eqref{APRI_SSOL_1}, we obtain
\begin{align}
\nonumber\left|\intO \psi'''(\varphi)|\grad\varphi|^2\delt\varphi\dx\right|&\leq C\left(1+\norml{\infty}{\varphi}^3\right)\normL{4}{\nabla\varphi}^2\norml{2}{\delt\varphi}\\
\nonumber &\leq C\left(1+\normh{3}{\varphi}^{\frac{3}{4}}\right)\normh{2}{\varphi}^{\frac{3}{2}}\norml{2}{\delt\varphi}\\
\nonumber &\leq C\left(1+\normh{3}{\varphi}\right)\left(\norml{2}{\varphi}+\norml{2}{\laplace\varphi}\right)\norml{2}{\delt\varphi}\\
\label{APRI_SSOL_13}&\leq C\left(1+\normh{3}{\varphi}^2\right)\left(1+\norml{2}{\laplace\varphi}^2\right)          +\frac{1}{4}\norml{2}{\delt\varphi}^2.
\end{align}
For the remaining term on the r.h.s. of \eqref{APRI_SSOL_9}, we observe that
\begin{equation*}
\chi\intO \nabla\sigma\cdot\nabla\delt\varphi\d x = \frac{\d}{\d t}\chi\intO\nabla\sigma\cdot\nabla\varphi\d x - \intO \nabla\delt\sigma\cdot\nabla\varphi \quad\text{for a.\,e. }t\in (0,T).
\end{equation*}
Collecting \eqref{APRI_SSOL_9}-\eqref{APRI_SSOL_11}, \eqref{APRI_SSOL_13} and using the last identity, we end up with
\begin{align*}
\frac{\d}{\dt}\frac{1}{2}\intO |\laplace\varphi|^2\dx  + \frac{1}{4}\intO |\delt\varphi|^2 &\leq C\left(1+\norml{2}{\divergence(\varphi\v)}^2  + \norml{2}{\Gamma_{\varphi}}^2 +\normh{3}{\varphi}^2\right)+C\left(1 + \normh{3}{\varphi}^2\right)\norml{2}{\laplace\varphi}^2\\
&\quad + \frac{\d}{\d t}\chi\intO\nabla\sigma\cdot\nabla\varphi\d x - \intO \nabla\delt\sigma\cdot\nabla\varphi\d x.
\end{align*}
Integrating this inequality in time from $0$ to $s\in (0,T]$, we obtain
\begin{align}
\nonumber\frac{1}{2}\norml{2}{\laplace\varphi(s)}^2 + \frac{1}{4}\norm{\delt\varphi}_{L^2(0,s;L^2)}^2 &\leq \norml{2}{\laplace\varphi_0}^2 + \int_{0}^{s}\alpha_1(t) + \alpha_2(t)\norml{2}{\laplace\varphi(t)}^2\d t\\
\nonumber&\quad + \chi\intO\nabla\sigma(s)\cdot\nabla\varphi(s)\d x - \chi\intO\nabla\sigma(0)\cdot\nabla\varphi_0\d x\\
\label{APRI_SSOL_15}&\quad  - \chi\int_{0}^{s}\intO \nabla\del_t\sigma\cdot\nabla\varphi\d x\d t,
\end{align}
where
\begin{equation*}
\alpha_1(t) \coloneqq C\left(1+\norml{2}{\divergence(\varphi\v)}^2 + \norml{2}{\Gamma_{\varphi}}^2  + \normh{3}{\varphi}^2\right),\qquad
\alpha_2(t)\coloneqq C\left(1+ \normh{3}{\varphi}^2\right).
\end{equation*}
Now, using \eqref{APRI_SSOL_1}, \eqref{APRI_SSOL_8} and $\varphi_0\in H_N^2$, we obtain
\begin{align*}
\left|\chi\intO\nabla\sigma(0)\cdot\nabla\varphi_0\d x \right| & = \left|\chi\intO \sigma(0)\laplace\varphi_0\d x\right|\leq C\left(1+\norml{2}{\Delta\varphi_0}^2\right),\\
\left|\chi\intO \nabla\sigma(s)\cdot\nabla\varphi(s)\d x\right| &\leq \norm{\sigma}_{C^0(H^1)}\sup_{s\in (0,T]}\norm{\nabla\varphi(s)}\leq C,\\
\left|\chi\int_{0}^{s}\intO \nabla\del_t\sigma\cdot\nabla\varphi\d x\d t\right|&\leq C\int_{0}^{s}\normh{1}{\del_t\sigma(t)}\d t.
\end{align*}
Together with \eqref{APRI_SSOL_15}, this implies
\begin{align}
\label{APRI_SSOL_16}\frac{1}{2}\norml{2}{\laplace\varphi(s)}^2 + \frac{1}{4}\norm{\delt\varphi}_{L^2(0,s;L^2)}^2 &\leq C\left(1+\norml{2}{\laplace\varphi_0}^2\right) + \int_{0}^{s}\beta_1(t) + \beta_2(t)\norml{2}{\laplace\varphi(t)}^2\d t,
\end{align}
where
\begin{equation*}
\beta_1(t) \coloneqq C\left(1+\normh{1}{\del_t\sigma}+\norml{2}{\divergence(\varphi\v)}^2 + \norml{2}{\Gamma_{\varphi}}^2 +  \normh{3}{\varphi}^2\right),\qquad\beta_2(t)\coloneqq C\left(1 + \normh{3}{\varphi}^2\right).
\end{equation*}
Due to \eqref{APRI_SSOL_1}, \eqref{APRI_SSOL_4} and \eqref{APRI_SSOL_8}, we know that $\beta_1,\beta_2\in L^1(0,T)$. Together with the assumption $\varphi_0\in H_N^2$, an application of Gronwall's lemma in \eqref{APRI_SSOL_16} yields
\begin{equation}
\label{APRI_SSOL_17}\norm{\laplace\varphi}_{L^{\infty}(L^2)} + \norm{\delt\varphi}_{L^2(L^2)}\leq C.
\end{equation}

\paragraph{Step 3:}
Combining \eqref{APRI_SSOL_1} and \eqref{APRI_SSOL_17}, from elliptic regularity theory we obtain
\begin{align}
\nonumber&\norm{\varphi}_{H^1(L^2)\cap L^{\infty}(H^2)\cap L^2(H^3) } + \norm{\sigma}_{H^1(H^1)\cap C^0(H^1)\cap L^{\infty}(H^2)} + \norm{\mu}_{L^2(H^1)\cap L^4(L^2)}\\
\label{APRI_SSOL_18}&\quad + \norm{\divergence(\varphi\v)}_{L^2(L^2)}  + \norm{\v}_{L^{\frac{8}{3}}(\H^1)}+ \norm{p}_{L^2(L^2)}\leq C.
\end{align}
Using \eqref{APRI_SSOL_18} and applying elliptic regularity theory in \eqref{WFORM_1b}, we obtain
\begin{equation}
\label{APRI_SSOL_19}\norm{\mu}_{L^2(H^2)}\leq C.
\end{equation}
Using a comparison argument in \eqref{WFORM_1c} and \eqref{APRI_SSOL_18}-\eqref{APRI_SSOL_19}, it follows that
\begin{equation}
\label{APRI_SSOL_20}\norm{\mu}_{L^{\infty}(L^2)}\leq C.
\end{equation}

\paragraph{Step 4:}
We now want to apply Lemma \ref{LEM_STOKES} with $q=2$. Using the generalised chain rule for Sobolev functions, \eqref{APRI_SSOL_18} and the assumptions on $\Gamma_{\v}$, it is straightforward to check that
\begin{equation}
\label{APRI_SSOL_21}\norm{\Gamma_{\v}}_{ L^{\infty}(H^1)}\leq C.
\end{equation}
Furthermore, since $\grad\varphi\in L^4(\L^{\infty}),~\mu+\chi\sigma\in L^{\infty}(L^2)$ with bounded norm, we observe that
\begin{equation}
\label{APRI_SSOL_22}\norm{(\mu+\chi\sigma)\grad\varphi}_{L^4(\L^2)}\leq C.
\end{equation}
Hence, using the assumptions on $\lambda(\cdot),~\eta(\cdot)$, an application of \eqref{STOKES_EST} yields
\begin{equation*}
\normH{2}{\v} + \normh{1}{p}\leq C\big(\eta_0,\eta_1,\lambda_0,\norm{\varphi}_{W^{1,4}}\big)\big(\normL{2}{(\mu+\chi\sigma)\grad\varphi} + \normh{1}{\Gamma_{\v}}\big) .
\end{equation*}
Integrating this inequality in time from $0$ to $T$, using \eqref{APRI_SSOL_21}-\eqref{APRI_SSOL_22} and recalling $\varphi\in L^{\infty}(W^{1,4})$ (since $H^2\hookrightarrow W^{1,4}$), we obtain
\begin{equation}
\label{APRI_SSOL_23}\norm{\v}_{L^4(\H^2)} + \norm{p}_{L^4(H^1)}\leq C.
\end{equation}
\paragraph{Step 5:}
Finally, due to the compact embedding $H^2\hookrightarrow C^0(\bar{\Omega})$ and because of \eqref{APRI_SSOL_18}, we obtain
\begin{equation*}
\norm{\varphi}_{C^0(\bar{Q})}\leq C.
\end{equation*}
For completeness, we summarize all the estimates we deduced in this Section, given by
\begin{align}
\nonumber&\norm{\varphi}_{ H^1(L^2)\cap C^0(\bar{Q})\cap L^{\infty}(H^2)\cap L^2(H^3)} + \norm{\sigma}_{H^1(H^1)\cap C^0(H^1)\cap L^{\infty}(H^2)} + \norm{\mu}_{L^{\infty}(L^2)\cap L^2(H^2)}\\
\label{APRI_SSOL_24}&\quad + \norm{\divergence(\varphi\v)}_{L^2(L^2)}  + \norm{\v}_{L^4(\H^2)}+ \norm{p}_{L^4(H^1)}\leq C.
\end{align}
These a-priori-estimates together with a Galerkin-scheme are enough to pass to the limit in the weak formulation to show existence of strong solutions. For the details when passing to the limit, we again refer to \cite{EbenbeckGarcke}, \cite{GarckeLam1}.
\paragraph{Step 6:} Since \eqref{MEQ_4} holds a.\,e. in $Q$, we see that $\varphi$ is a solution of
\begin{alignat*}{3}
\laplace\varphi &= \psi'(\varphi) - \mu - \chi\sigma&&\quad\text{a.\,e. in }Q,\\
\nabla\varphi\cdot\n &= 0&&\quad\text{a.\,e. on }\Sigma.
\end{alignat*}
Since $\psi'(\varphi) - \mu - \chi\sigma\in L^2(H^2)$ with bounded norm, elliptic regularity theory implies
\begin{equation}
\label{APRI_SSOL_25}\norm{\varphi}_{L^2(H^4)}\leq C.
\end{equation}
Due to the continuous embedding $L^{\infty}(H^1)\cap L^2(H^3)\hookrightarrow L^8(L^{\infty})$ and by \eqref{APRI_SSOL_24}, this implies $(\mu+\chi\sigma)\nabla\varphi\in L^8(\L^2)$ with bounded norm. Consequently, with the same arguments as used for \eqref{APRI_SSOL_23}, we deduce that
\begin{equation}
\label{APRI_SSOL_26}\norm{\v}_{L^8(\H^2)} + \norm{p}_{L^8(H^1)}\leq C,
\end{equation}
which completes the proof.
\section*{Acknowledgements} 
This work was supported by the RTG 2339 \grqq Interfaces, Complex Structures, and Singular Limits" of the German Science Foundation (DFG). The support is gratefully acknowledged. 
\printbibliography[heading=bibintoc,title={References}]
\end{document}